\documentclass[12pt]{article}
\usepackage {amsfonts, amsthm,euscript, amscd, amsmath,amssymb, color}
\usepackage{epsfig}
\usepackage{graphicx}
\author { H. Abels, G.A. Margulis and G.A. Soifer}
\title {The Auslander conjecture for dimension less than 7 (RF) }
\normalfont

\def\nl{\newline}

\def\ltextindent#1{\hbox to \hangindent{#1\hss}\ignorespaces}
\def\({\left (}
\def\[{\left [}
\def\){\right )}
\def\]{\right ]}

\def\br{\Bbb R}

\def\g{\gamma}
\def\G{\Gamma}

\def\o{\omega}

\def\ve{\varepsilon}

\def\Aff{\text{Aff\,}}

\def\la{\longrightarrow}

 \font\normalfont=cmr12 at
12truept   at 20.74truept

\def\la{\longrightarrow}

\def\pro#1{\medskip{\bf #1 }\it}
\def\endpro{\normalfont\medskip}

\begin{document}
\maketitle
\noindent \textbf {Abstract.}
In 1964  L. Auslander
conjectured that every crystallographic subgroup $\Gamma$ of an affine group $\text{Aff}(\mathbb R^n)$
is virtually solvable, i.e. contains a solvable subgroup of finite index. D.~Fried and W.~Goldman proved Auslander's conjecture for
$n=3$ using cohomological arguments.  We prove the Auslander conjecture for $n < 7$. The proof is based mainly on dynamical arguments.
In some cases we use the cohomological argument which we could avoid but it would significantly lengthen the proof.  

\section{Preamble}\label{0} Our paper with the same title was put on the archive in 2012[AMS5]. Motivated by questions of readers we wrote a new version. Here we explain ideas of proofs in more detail. It makes the text longer 
but we hope that this version is more reader friendly (RF). It is submitted for publication. We herewith make it available to the public.

\section{Introduction}\label{0} 

Let $X$ be a topological space and  $\G$ be a subgroup of the group of homeomorphisms of $X$. A subgroup $\Gamma$ is said to act  {\it
properly discontinuously} on $X$ if for every compact
subset $K$ of $X$ the set $\{g\in\Gamma:  g K\cap
K\ne\emptyset\}$ is finite.  A subgroup $\Gamma $  is called {\it crystallographic} if $\Gamma$
acts properly discontinuously on $X$ and the orbit space
 $\Gamma \setminus X$  is compact.  
  \par  The study of crystallographic groups has a long history.
The crystallographic groups of hyperbolic transformations have been investigated by Fricke and Klein in the lectures on the theory of automorphic functions [FK]. In $3$-dimensional Euclidean space Fedorov [F],  Schoenflies [Sc] and later Rohn [Ro] have shown that there are only a finite number of essentially different kinds of euclidean crystallographic groups. The $3$-dimensional euclidean crystallographic groups are the symmetry groups of crystalline structures and so are of fundamental importance in the science of crystallography. 
\par D. Hilbert wrote in his famous lecture delivered on the IMC at Paris, 1900 ([Hil], 18th Problem)
:\\
" Now, while the results and methods of proof applicable to elliptic and hyperbolic space hold directly for n-dimensional space also, the generalization of the theorem for euclidean space seems to offer decided difficulties. The investigation of the following question is therefore desirable: \\
\textit{Is there in n-dimensional euclidean space also only a finite number of essentially different kinds of groups of motions with a fundamental region? "}
\par In response to this problem  Bieberbach  showed in  a series  of  papers   that this was so. The key result is the following  famous theorem of Bieberbach. A group $\G$ acting isometrically on  the $n$--dimensional Euclidean space $\mathbb R^n$ with compact
quotient  $\Gamma\setminus\mathbb{R}^n$ contains a subgroup  of
finite index consisting of translations. In particular, such a
group $\Gamma$ is virtually abelian, i.e. $\Gamma$ contains an
abelian subgroup of finite index. Moreover, in  [B1, B2, B3] Bieberbach proved that a group $\G$ acting isometrically and properly discontinuously on  the $n$--dimensional Euclidean space is virtually abelian. 
\par A natural way to generalize the classical problem is to broaden the class of allowed motions and consider crystallographic groups of affine transformations. This raises the question of the group-theoretic properties of affine crystallographic groups.
\par Let us recall that the group  $G_n = \text{Aff}(\mathbb R^n)$ of affine transformations of
$\mathbb R^n$ is the semidirect product $GL_n(\mathbb R)\ltimes
\mathbb R^n$ where $\mathbb R^n$ is identified with the  group of
its translations.  Let $l:G_n \rightarrow GL_n(\mathbb{R})$ be the natural
homomorphism. Then  $l(g)$ is called the \textit{linear part} of the affine transformation $g$. Let $X \subseteq G_n $, then the set $l(X)= \{l(x), x \in X\}$ is called the \textit {linear part} of $X$.\\
Auslander proposed  the following conjecture in [Au]. \\
\pro {\it The Auslander Conjecture.} Every crystallographic subgroup $\Gamma$ of $G_n$
is virtually solvable, i.e. contains a solvable subgroup of finite
index. \endpro\\
The proof in [Au] of this conjecture is unfortunately incorrect, but the conjecture is still an open and central problem (see Milnor [Mi2]). 
\par It is easy to see that there exists a nilpotent, non virtually abelian affine crystallographic group. It is well known [Mo] that every discrete virtually solvable linear  group and in particular every virtually solvable  discrete subgroup of $G_n$  is virtually polycyclic.  J. Milnor showed that every virtually polycyclic group can act properly discontinuously and affinely on some affine space [Mi1]. 
\par There is an additional geometric interest in properly discontinuous  groups
since they can be represented as fundamental groups of manifolds
with certain geometric structures, namely complete flat affine manifolds. If
$M$ is a complete flat affine manifold, its universal covering manifold is
isomorphic to $\mathbb{R}^n$. It follows that its fundamental 
group $\Gamma= \pi_1(M)$ is in a natural way a properly
discontinuous  torsion-free subgroup of $G_n$.
Conversely, if $\Gamma$ is a properly discontinuous
torsion-free subgroup of $G_n$, then $\Gamma
\setminus\mathbb{R}^n$ is a complete flat affine manifold $M$ with
$\pi_1(M) = \Gamma$. Therefore every virtually polycyclic group is the fundamental group of a complete affinely flat manifold by Milnor's theorem mentioned above.  J. Milnor proposed the following questions in [Mi1]: \bigskip\\
\pro {\it Question 1} Let $\G$ be a torsion free virtually polycyclic group of rank $k$. Does there exist a $k$--dimensional compact complete affinely flat manifold $M$ with $\pi_1(M) \cong \G ?$ 
\endpro\\
Now it is known that not every  finitely generated nilpotent group is an  affine crystallographic group  [B]. This gives a negative answer to Question 1. \bigskip\\
\pro {\it Question 2}  Does there exist a complete affinely flat manifold $M$ such that $\pi_1(M) $ contains a free group ? \endpro \\
In comments to the second question Milnor wrote: "I do not know if such a manifold exists even in dimension 3" and proposed "to construct a Lorentz-flat  example by starting with  a discrete subgroup  $\mathbb{Z}\ast \mathbb{Z} \leq SO(2,1)$  then  adding  translation  components to  obtain  a  group  of  isometries of  Lorentz  3-space; but  it  seems difficult  to  decide whether  the  resulting group action is properly  discontinuous"  [Mi2, p. 184].  
\par G. Margulis gave a positive answer to Question 2 in dimension 3 in [M]. He constructed a free non-abelian subgroup $\G$ of  isometries of  Lorentz  3-space acting properly discontinuously on 
$\mathbb{R}^3.$  Clearly $l(\G ) \subseteq SO(2,1).$ Then  we proved in [AMS3, Theorem B] that for a non degenerate form $B$ of signature $(p, q)$ where $p=q+1$ and $q$ is odd, there exists a free group $\G \leq $ \text{Aff} $(\mathbb{R}^{2q+1})$ acting properly discontinuously such that the linear part $l(\G) $ of $\G$ is Zariski dense in $SO(B)$. Therefore in any dimension $n \geq 3$ there exists a complete affinely flat manifold $M$ such that $\pi_1(M) $ contains a free non- abelian group. 
\par Let $B$ be a non degenerate quadratic form. Set  $G_B= \{x \in G_n: l(x) \in O(B) \}.$ Clearly $G_B = O(B) \ltimes
\mathbb R^n.$ Let $\G$ be an affine crystallographic group, and suppose $\G \subset G_B$ where $B$ has signature $(p,q)$.  Remark, that if  $q =0$ we have the case of isometric affine actions. D. Fried and W. Goldman in [FG]  proved, that if $\G$ is a crystallographic subgroup of  $G_B$ where $B$  is a non degenerate quadratic form of signature $(2,1)$ then
$l(\G)$ is not Zariski dense in $O(2,1).$  They use this theorem to deduce the Auslander conjecture for dimension at most 3.   W.Goldman and Y. Kamishima  proved in [GK] that a crystallographic subgroup of  $G_B$ is virtually solvable for $q=1$.  In [AMS4] we proved that a crystallographic group $\G \subseteq G_B$ is virtually solvable if $B$ is a quadratic form of signature $(p,2).$ F. Grunewald and G. Margulis  proved in [GM] that if the linear part is a subgroup of a simple Lie group of real rank 1, then $\G$ is virtually solvable. This result was generalized in  [To1]. Namely it was proved that if $l(\G) \subseteq G$ and the semisimple part of $G$ is a simple group of real rank 1, then $\G$ is  virtually solvable. Finally in [S2] and [To2] it was proved, that if  $l(\G) \subseteq G$ and  every non-abelian simple subgroup of $G$ has real rank $\leq 1$  then $\G$ is virtually solvable. Let us remark, that all papers [FG], [GK], [GM], [S2] and [To1,2]  basically use the same idea which was first introduced in [FG]. We call this idea "the cohomological argument" because it is based on using the virtual  cohomological dimension of $\G.$  By contrast,  [AMS4] and [M] are based on a completely different approach, namely on dynamical ideas (see also [AMS1,2,3]). 
\par  In [To3] the author attempts to prove the Auslander conjecture for  dimensions 4 and 5.  Unfortunately, the proof there is incomplete. 
Thus the only dimensions for which there is a complete proof of the Auslander conjecture for $G_n$  are $n \leq 3$ (see  [FG, Theorem, section 2.13] ) 
\par Let us mention the following result due to M. Gromov [Gr].
Let $M$ be a connected compact Riemannian manifold. Denote by $d=d(M)$
the diameter of $M$, and by $c^+ = c^+(M)$ and $c^- = c^-(M)$, respectively,
the upper and lower bounds of the sectional curvature of $M$. We set
$c =c(M)= \max ( c^+, c^- ).$ We say that $M$ is almost $ \varepsilon $--flat,
$\varepsilon \geq 0$ if $cd^2 < \varepsilon.$ Then for sufficiently small $\varepsilon$ the
fundamental group of an $\varepsilon$--flat manifold is a virtually nilpotent group, i.e. contains
a nilpotent subgroup of finite index. This result again shows that $M$ being close to Euclidean has strong implications for the algebraic structure
of the fundamental group $\pi_1 (M).$
\par In [DG1], [DG2], [CDGM] and [Me] there were studied properly discontinuous subgroups $\Gamma$ of
the affine group \text{Aff}$(\mathbb{R}^3)$   whose linear part $l (\Gamma)$  leaves a quadratic form of signature $(2,1)$ invariant. 
\par The aim of this paper is to prove  the following theorem which was announced in [AMS2]\\

\pro{\it Main Theorem} Let $\G$ be a crystallographic subgroup of \text{Aff}$(\mathbb{R}^n)$ and $n < 7$, then $\G$ is virtually solvable.
\endpro\\
The proof of this theorem is based mainly on dynamical arguments. In some cases we use the cohomological argument to shorten the proofs.
\par Let us give a short description of the paper. In section 2 we introduce the terminology we will use throughout the  paper and recall some basic results about the dynamics of the action of hyperbolic elements. We show that every element of the connected component of the Zariski closure of an affine group acting properly discontinuously has one as an eigenvalue. This simple but useful fact will be used in section 3. As the first step in section 3 we obtain a list of all possible semisimple groups $S$ which might be a semisimple part of the Zariski closure of an affine group $\G $ that acts properly discontinuously for $ n \leq 6$ and does not have $SO(2,1)$ as a quotient group.  Using this list we prove the Auslander conjecture in dimension 4 and 5 in section 4. Actually we show a bit more. Namely, if the semisimple part of the Zariski closure of $\G$ is one from the list, then $\G$ does not act properly discontinuously. Note that based on [M]  it is not difficult to show,  that if  the semisimple group $S$ contains $SO(2,1)$ as a quotient group, there exists an affine group $\Gamma$ acting properly discontinuously such that the linear part of $\Gamma$ is Zariski dense in $S$. 
In section 5 we show that the  semisimple part $S$ of the Zariski closure of $l(\G)$ cannot be $SO(3,2)$ or $SO(3) \times SL_3(\mathbb{R}).$ The proof is based on the cohomological argument we have mentioned above. Namely, we will compare the virtual cohomological dimension of $\G$ and the dimension of the symmetric space $S/K$, where $K$ is a maximal compact subgroup of $S.$ We will prove that none of these cases is possible.
\par The most difficult part is to show that the semisimple part of the Zariski closure of $l(\G)$ is not $SO(2,1) \times SL_3(\mathbb{R}).$. This is done in section 6. We show that  it is possible to change the sign of a hyperbolic element (see  Main Lemma 6.7) in this case.  Thus, by Lemma 6.5, we conclude that the semisimple part of the Zariski closure of $l(\G)$ cannot be  $SO(2,1) \times SL_3(\mathbb{R}).$  Hence none of the possible non-trivial semisimple groups can be the semisimple part of the Zariski closure of $\G$. Therefore the semisimple part of the Zariski closure of $\G$ is trivial. Hence $\G$ is virtually solvable. 
\par In the final section 7 we discuss Auslander's conjecture in dimension 7 and state two open problems. We believe that  answers (positive or negative) to these questions are 
essential for further progress on Auslander's conjecture.\\
The authors would like to thank several institutions and foundations for their
support during the preparation of this paper: Bielefeld
University, the Emmy Noether Institute, Bar-Ilan University,  Yale University,
SFB 701 "Spektrale Strukturen und Topologische Methoden in der Mathematik", National Scientific Foundation under grant DMS-0801195 and
 DMS-1265695, USA-Israel Binational Science foundation under grants 2004010 and 2010295,
Israel Science foundation under grant  657/09, Max Planck Institute for Mathematics (Bonn).
Without all these supports, the paper whose authors live on three different continents could not have seen the light of day. 
\section{Dynamical properties of the  action of hyperbolic elements}\label{0}
\pro {\it 2.1. Notation and terminology.}\endpro In this section
we introduce the terminology we will use throughout the 
paper. We also prove and recall some basic results about the dynamics of the action of
hyperbolic elements [A], [AMS1, 4]. \\
 An affine transformation $x \in G_n$ is a pair of a linear transformation $l(x)$ and a vector
$v_x \in \mathbb{R}^n$. Denote by $0$ the point of origin. Let $x0= 0 + v_x.$  By definition for every point $q$ of the affine space $\mathbb{R}^n$ we have $x q = l(x) \overrightarrow{0q} +v_x+ 0.$  Let $x$ be an affine transformation such that $v_x= 0.$ Hence $x0 = l(x)\overrightarrow{0}+0 =0.$ Therefore in the case $v_x= 0$ we will sometimes write $l(x)q$ instead of $xq$ for any $q \in \mathbb{R}^n.$\\
There exists a linear representation  $\phi : G_n  \rightarrow GL_{n+1}(\mathbb{R})$ defined by
$$ \phi : x \mapsto 
\left(
\begin{array}{cc}
l(x)& v_x\\
0&1\\
\end{array}
\right)\,
 \eqno (\ast)$$ 
 Hence we can and will consider $G_n$ as a linear group. \\
 Let $\Gamma$ be an affine group and let $G$  be the Zariski closure of $\Gamma.$  Denote by $\tilde{S}$ a reductive part  of $G.$ There exists an element $t \in G$ such that $t \tilde{S} t^{-1} \leq l(G_n).$ Set $q_0 = t^{-1} 0.$ Then $\tilde{S} q_0 =q_0.$ \\ 
\noindent \textbf{2.2.} Let $V$ be a finite dimensional vector space over a local field $k$ of characteristic $0$ with absolute value  $| \cdot|.$
Let $g \in GL(V)$ be a linear transformation and let $
\chi_g(\lambda)= \prod^n_{i=1}(\lambda - {\lambda}_i) \in
k[\lambda] $ be the characteristic polynomial of $g$. For arbitrary positive $\alpha \in \mathbb{R}^+$ we set $\Omega_{\alpha}(g)= \{ \lambda_i : |\lambda_i| > \alpha
 \}.$ Put $\chi_{\alpha} (\lambda)= \prod_ {\lambda_i \in \Omega_{\alpha}(g)}(\lambda -\lambda_i).$
Then $\chi_{\alpha} (\lambda)$ belong to $k[\lambda]$ since the absolute value of an
 element is invariant under Galois automorphisms. Therefore $\chi_{\alpha} (g) \in GL(V).$ 
 Set 
  $\overline{\Omega}_1(g)= \{ \lambda_i, 1 \leq i \leq n : |\lambda_i| \geq 1 \} .$ 
 Let $\overline {\chi}_1 (\lambda)= \prod_ {\lambda_i \in \overline{\Omega}_{1}(g)} (\lambda -\lambda_i).$ By the same arguments as above we conclude that  $\overline {\chi} _1(\lambda)  \in k[\lambda].$ Therefore
$ \overline {\chi}_1 (g) \in GL(V). $   Put  $A^+(g) = \ker \chi_{1} (g)$,  $D^+(g) = \ker {\overline{\chi}}_{1} (g),$
$A^-(g) =A^+(g^{-1}),$ $D^-(g) =D^+(g^{-1})$ and $A^0(g) =D^+(g) \cap D^-(g).$  
\par  We set $A^+(g) =A^+(l(g)),$  $D^+(g) =D^+(l(g)),$ $A^-(g) =A^-(l(g)),$  $D^-(g) =D^-(l(g))$ and $A^0(g) =A^0(l(g))$ for an affine transformation  $g \in G_n$ [AMS4].
\noindent\textbf{2.3.}  Let $g$ be an element in
$GL({\mathbb{R}}^n)$. Then the space $V=\mathbb R^n$ is the direct sum of the three $g$-invariant subspaces $A^+(g),
~{A^-(g)}$ and ${A^0(g)}.$ By definition, all
eigenvalues of the restriction $g|_ {A^+(g)}$(resp.~$g|
_{A^-(g)}$, ~$g|_{A^0(g)}$) have absolute value greater than 1
(resp. less than 1, equal to 1).  Clearly $D^+(g)= A^+(g) \oplus
{A^0(g)}$ and $D^-(g)= A^-(g) \oplus {A^0(g)}$. 
\par Let $G$ be a subgroup of $GL(V)$ and let $g \in G.$ Set $V^0_g = \{v \in V ; gv=v\}.$ A semisimple element $g \in G$ is called \textbf {regular }in $G$ if $\dim V_g^0= \min \{\dim V^0_t| $ $\, t \in G, t$ \text {is a semisimple element}$ \}$.
\par If for a semisimple element $g \in G$ we have $\dim(A^0(g))= \min \{\dim A^0(t)|$ $ \, t \in G, t$ \text {is a semisimple element}$ \}$, then $g \in G$ is called $\mathbb{R}$-\textbf{regular} in $G \leq GL(V).$ Let $G$ be an affine group, $G < G_n.$  An affine transformation $g \in G$ is called regular  (respectively $\mathbb{R}$-regular) if $l(g)$ is a regular (respectively $\mathbb{R}$-regular) element of $l(G).$ \par Our definition of $\mathbb{R}$-regular element slightly differs from that of [P] were it was first introduced. Note that the set of $\mathbb{R}$-regular elements in an algebraic group $G$ need not be Zariski open in $G.$ Nevertheless under some conditions a Zariski dense subgroup of 
 an algebraic group $G$ contains an $\mathbb{R}$-regular element 
 [P],[AMS1],[AMS4]. For example this is true if $G=SO(B)$ where $B$ is a non degenerate form of signature $(p,q)$ and $\Gamma$ is a Zariski dense subgroup of $G$. Note that in case $p=2,q=1$ every hyperbolic element is $\mathbb{R}$-regular.
 \par 
 \noindent \textbf{2.4.} If we use topological concepts and do not specify the topology, we mean the Zariski topology. If we refer to the Euclidean topology we mention this explicitly and use expressions like Euclidean-open, Euclidean-connected, etc.\\
Let  $\|\centerdot \|$ and $d$ denote the norm and metric on
$\mathbb {R}^n$ corresponding to the standard inner product on
$\mathbb {R}^n$. Let $\|g\|_-$ be the norm of the restriction $g|_{A^-(g)}.$ Denote by $\|g\|_+ =\|g^{-1}\|_-$ and put $s(g)= \max \{\|g\|_+, \|g\|_-\}.$ A regular element $g$ is called $\it{hyperbolic}$ if $s(g) < 1.$ It is clear that for a regular element $g$ there exists a number $N$ such that for $n > N$ the element $g^n$ is hyperbolic.\par Let $P= \mathbb{P}(\mathbb{R}^n)$ be the projective space corresponding to $\mathbb {R}^n$.
Let $\pi : \mathbb {R}^n \setminus \{0\} \rightarrow P$ be the natural projection. For a subset $X$ of $\mathbb{R}^n$  we denote $\pi(X)=\pi (X \setminus\{0 \})$ .
\par The metric $\|\centerdot \|$ on  $\mathbb{R}^n$ induces the standard metric $\widehat{d}$ on the projective space $P = \mathbb{P}(\mathbb{R}^n)$ by the formula (see [T])
$$\widehat{d}(p,q) = \frac{\|v \wedge w \|}{\|v \|\cdot\|w\|}, p =\pi (v), q=\pi(w)$$. Thus for any point $p \in P$ and any subset $A \subseteq P,$ we can define
$$\widehat{d}(p,A) = \inf_{a\in A} \widehat{d}(p,a).$$
Let $A_1$ and $A_2$ be two subsets of $P$. We define
$$\widehat{d}(A_1,A_2) = \inf_{a_1\in A_1}\inf_{ a_2 \in A_2} \widehat{d}(a_1,a_2) $$ and
$$\widehat{\rho}(A_1, A_2) = \inf \{ R; A_2 \subseteq B( A_1, R),  A_1 \subseteq  B(A_2, R)\}$$
where $B(A, R) =\bigcup _{a\in A}  B (a, R) .$\\
For two subspaces $W_1 \subseteq \mathbb{R}^n $ and $W_2 \subseteq \mathbb{R}^n$ we put $\widehat{d}(W_1,W_2)= \widehat{d}(\pi(W_1 ), \pi(W_2))$ and  $\widehat{\rho}(W_1,W_2)= \widehat{\rho}(\pi(W_1), \pi(W_2))
.$ A hyperbolic element $g $ is called $\varepsilon$-\textit{hyperbolic} if  $$\widehat{d}(A^+(g), D^-(g)) \geq \varepsilon$$ and  $$\widehat{d}(A^-(g),D^+(g)) \geq \varepsilon.$$  Two different hyperbolic elements $g_1$ and $g_2$ are called \textit{transversal} if\\ $A^{\pm}(g_1)\bigcap  D^{\mp}(g_2) = \{0\}$ and $A^{\pm}(g_2)\bigcap  D^{\mp}(g_1) = \{0\}$. Let $B$ be a  non degenerate quadratic form and let $g \in SO(B)$ be a $\mathbb{R}$-regular element. Since $A^+(g)$ (resp.$A^-(g)$ ) is the unique maximal
isotropic subspace of $D^+(g)$ (resp. $D^-(g)$) it is easy to see that two hyperbolic $\mathbb{R}$-regular elements $g_1$ and $g_2$ of $SO(B)$ are transversal if and only if $A^+(g_1) \bigcap A^-(g_2) =\{0\}$ and
$A^+(g_2) \bigcap A^-(g_1) =\{0\}.$
\par Clearly $g$ and $g^{-1}$ are not transversal for any regular element $g$. Nevertheless it is quite important to be able to find an element $t $ of a given linear group $G$ such that $g$ and $t g^{-1} t^{-1}$ are transversal. It is possible for example for $G =SO(B).$  For an arbitrary $G < GL(V)$ it might happen that  $\dim A^+(g) \neq \dim  A^-(g).$ Assume for example that  $\dim A^+(g)  > \dim  A^-(g) .$ Thus for every $t \in G$ we have $tA^+(g) \cap D^+(g) \neq \varnothing.$ Therefore there is no $t \in G$ such that  $g$ and $t g^{-1} t^{-1}$ are transversal. This is the motivation of the following definition. \\
\noindent\pro {\it Definition 2.5} Let $g \in G \subseteq GL(V)$ be 
a regular element  such that $\dim A^+(g) \geq \dim A^-(g).$ We will say that $g$ can be transformed into a transversal pair inside $G$ if there exists an element $t \in G$ and a 
 subspace $W \subset A^+(g)$ such that 
$V =t W  \oplus D^+(g).$ \endpro \\
\noindent\pro {\it Remark.}\endpro It is easy to see that an element $g \in G$ can be transformed into a transversal pair inside $G$ if and only if there exists an element $t \in G$ such that $D^+(g)+ tA^+(g) =V.$
\par The next proposition shows that this property depends only on the Zariski closure $\overline{ G}$ of a group $G$, and thus $G$ can be safely ignored in most of what we do.
\\
\noindent\pro{\it Proposition 2.6} Let $\overline{ G}$ be the Zariski closure of $G \subseteq SL(V).$ Assume that $\overline{ G}$ is connected non-solvable group.  
Let $g \in G $ be a regular element of  $\overline{ G}$ which  can be transformed into a transversal pair inside  $\overline{ G}.$ Thus there exist a subspace $W$ of $A^+(g)$ and $t \in \overline{G}$ such that $V = D^+(g) \oplus tW .$ Then 
\begin{enumerate}
\item[(A1)] The set $\Omega(g) =  \left\lbrace  t \in \overline{G} , V = D^+(g) \oplus t W \right\rbrace  $ is non-empty and open,
 \item[(A2)] Let $\Omega_{ab} (g)$ be the set of all $ t \in \overline{G}$ such that $ g $ and $ tgt^{-1 }$  do not commute. 
 Then the set $\Omega_{ab} (g) $ is open. 
 \item[(A3)] If the set $\Omega_{ab} (g) $ is non-empty there exists $t \in G\cap \Omega(g) \cap \Omega_{ab} (g).$ Therefore we have  $V = D^+(g) \oplus tW $ and the group generated by $g$ and $tgt^{-1}$ is not commutative. 
 \item[(A4)] The set $$\Omega= \left\lbrace (t,g), t \in \overline{G}, g \in \overline{G}: tA^+(g) +D^+(g) =V, g\,\text{is a regular elment of}\,\overline{G} \right\rbrace $$ is non- empty and  open in $\overline{G} \times \overline{G}.$
\end{enumerate} 
 \endpro 
 \noindent\pro {\it Proof.} \endpro The sets $\Omega(g)$ and $\Omega_{ab} (g)$ are Zariski open since their complement is determined by algebraic equations. From Definition 2.5 follows that  $\Omega(g) \neq \varnothing.$ The semisimple part of $\overline{ G}$ is not trivial, therefore the set $\Omega_{ab}(g) \neq \varnothing.$ This proves (A1) and (A2). Clearly $\Omega$ is the intersection of two open subsets of $\overline{G}\times \overline{G}.$ Thus $\Omega$ is an open subset. Since there exists a regular element of  $\overline{ G}$ which  can be transformed into a transversal pair inside  $\overline{ G}$ we conclude that $\Omega \neq \varnothing.$ Note that $G$ is dense and $\Omega(g)$ and $\Omega_{ab}(g)$ are open subsets in $\overline{ G}.$ Hence the set $\Omega(g) \cap \Omega_{ab}(g) \cap G$ is non-empty. This proves the proposition.
\par Two transversal hyperbolic elements $g_1$ and $g_2$ of $GL(V)$ are called $\varepsilon$- \textit{transversal}, $$ \min_{1 \leq i \neq j \leq 2}\{\widehat{d}(A^+(g_i), D^-(g_j)),\widehat{d}(A^-(g_i), D^+(g_j))\} \geq \varepsilon .$$
 Let $g_1 \in SO(B)$ and $g_2 \in SO(B)$ be two hyperbolic elements $g_1$ and $g_2.$ As we pointed out two elements $g_1 \in SO(B)$ and $g_2 \in SO(B)$ are transversal 
if and only if $A^+(g_1) \cap A^-(g_2) = \{0\}$ and $A^-(g_1) \cap A^+(g_2) = \{0\}.$  Thus $g_1$ and $g_2$ are transversal if and only if $g_1^{-1}$ and $g_2^{-1}$ are transversal. Then
\begin{enumerate}
\item[(1)] for every $\varepsilon$ there exists $\delta =\delta(\varepsilon)$ such that $g_1$ and $g_2$ are $\varepsilon$-transversal if $\widehat{d} (A^+(g_1), A^-(g_2)) > \delta$ and $\widehat{d}(A^-(g_1), A^+(g_2)) > \delta$.
\item[(2)] For every $\delta $ there exists $\varepsilon=\varepsilon(\delta)$ such that if   $\widehat{d} (A^+(g_1), A^-(g_2)) > \delta$ and $\widehat{d}(A^-(g_1), A^+(g_2)) > \delta$ 
then the two hyperbolic elements $g_1$ and $g_2$ are $\varepsilon$ -transversal.
\end{enumerate}
\noindent Clearly   
\begin{enumerate}
\item[(3)]  $g_1$ and $g_2$ are $\varepsilon$-transversal  if and only if $g_1^{-1}$ and $g_2^{-1}$ are $\varepsilon$-transversal.
\end{enumerate}
\par An affine transformation is called $\textit {hyperbolic}$  (respectively $\varepsilon$-\textit{hyperbolic})\, if $l(g)$ is hyperbolic  (respectively $\varepsilon$-\textit{hyperbolic})).
Recall the following useful Lemma 
\bigskip
\nl
\noindent\pro{\it Lemma 2.7. [AMS3] } There exists $s(\ve)<1$ and $c(\ve)$
such that for any two $\ve$-hyperbolic $\ve$-transversal
elements $g, h \in GL(V)$  with $s(g)<s(\ve)$ and $s(h)<s(\ve),$ we have \\
{\rm (1)}~ the element ~$gh$ is $\ve/2$-hyperbolic and is
$\ve/2$-transversal to both $g$~and ~$h$;\nl
{\rm (2)}~ $\widehat{\rho}(A^+(gh), A^+(g)) \leq c(\ve)s(g))$;\nl
{\rm (3)}~$\widehat{\rho}(A^-(gh), A^-(h)) \leq c(\ve)s(h))$;\nl
 {\rm (4)}~$s(g h)\leq c(\ve)s(g)s(h).$
\endpro
\bigskip\\
\pro{Proposition 2.8} Let $\G$ be an affine group acting properly discontinuously.  Let $g $ be an element of  the connected component of the Zariski closure of $\G .$  Then $l(g)$  has $1$ as an eigenvalue.
\endpro\\
\noindent\pro {\it Proof} \endpro It is easy to see that for $x \in G_n$,  if  $l(x)$ does not have $1$ as an eigenvalue then x has a fixed point. Thus every element of an affine torsion free group acting properly discontinuously  has one as an eigenvalue. Let $G$ be the Zariski closer of $\Gamma$ and let $G^0$ be the connected component of $G$. It is well known, that exists a finitely generated subgroup $\G_0$ of $\G$  such that the Zariski closure of $\G$ and $\G_0$ coincide. Since $G^0$ is an open and closed subgroup of $G$ a finite index subgroup $\Gamma_1 = \Gamma_0 \cap G^0$ of $\Gamma_0$ is a finitely generated group which is dense in $G^0.$ By  Selberg's lemma we conclude that there exists a  torsion free subgroup $\G_2 \leq \G_0$ of finite index.Thus $\Gamma_2$ is Zariski dense in $G^0$.Since $\G_2$ acts properly discontinuously the linear part $l(x)$ has one as an eigenvalue for every $x \in \G_2$. Consequently the same is true for every element of the Zariski closure of $\G_2.$ This proves the statement.
\par This simple proposition  will help us to list all possible semisimple connected Lie groups which can be a semisimple part of the Zariski closure of a subgroup of $G_n, n  \leq 6$, acting properly discontinuously (see section 4). 
\par Let $\G$ be an affine group acting properly discontinuously. Obviously, $\G$ is a crystallographic group if and only if every finite index  subgroup of $\G$ is crystallographic. Thus, since the connected component of a Zariski closed group is a subgroup of finite index, we will  assume from now on that the Zariski closure of $\G$ is connected. Therefore every element  of $\G$ has one as an eigenvalue by Proposition 2.8. Hence for a semisimple element $g$ of $\G$ there exists a $g$-invariant line $L_g$. The restriction of $g$ to $L_g$ is
the translation by a non- zero vector $t_g$. Let us note that all
such lines are parallel and the vector $t_g$ does not depend on the
choice of $L_g$.  We take for $g$ the $g$-invariant line $L_g$  that is closest to the origin. 
Let us define the following affine subspaces: $E_g^+=D^+(g)+L_g$,
~$E_g^-=D^-(g)+L_g, ~E_g^+ \cap  E_g^-= C_g.$  Let $p \in L_g$ be a point. Then $t_g  = \overrightarrow{p\, gp} .$ Clearly $t_g =- t_{g^{-1}},$ $L_g=
L_{g^{-1}}.$ For $h =t g t^{-1}$ we have $$t_h =l(t)t_g.  \eqno (**)$$ 
\noindent\pro {\it Proposition 2.9} Let $G \subset GL(V) $ be the Zariski closure of the linear part of an affine group $\Gamma.$ Let $S$ be a semisimple part of $G$ and let $U$ be the unipotent radical of $G$. Assume that  $G$ is a connected group and $V$ is the direct sum of two non-trivial $S$-invariant subspaces $V_0$ and $V_1$ with the following properties. 
\begin{enumerate}
\item[(1)] $gv =v$ for all $g \in G, v \in V_0$ and the induced action $g : V/V_0 \rightarrow V/V_0 $ is trivial for all $g \in U.$
\item[(2)] The restriction $g|_{V_1}$ for one (then for every) regular element $g$ of $S$ does not have one as eigenvalue.
\item[(3)] Every regular element $s \in S$ can be transformed into a transversal pair inside $S$. 
\end{enumerate}
Then the group $\Gamma$ does not act properly discontinuously \endpro  \\
\begin{figure}[!h] 
\label{pic-gamma}
 \centerline{ 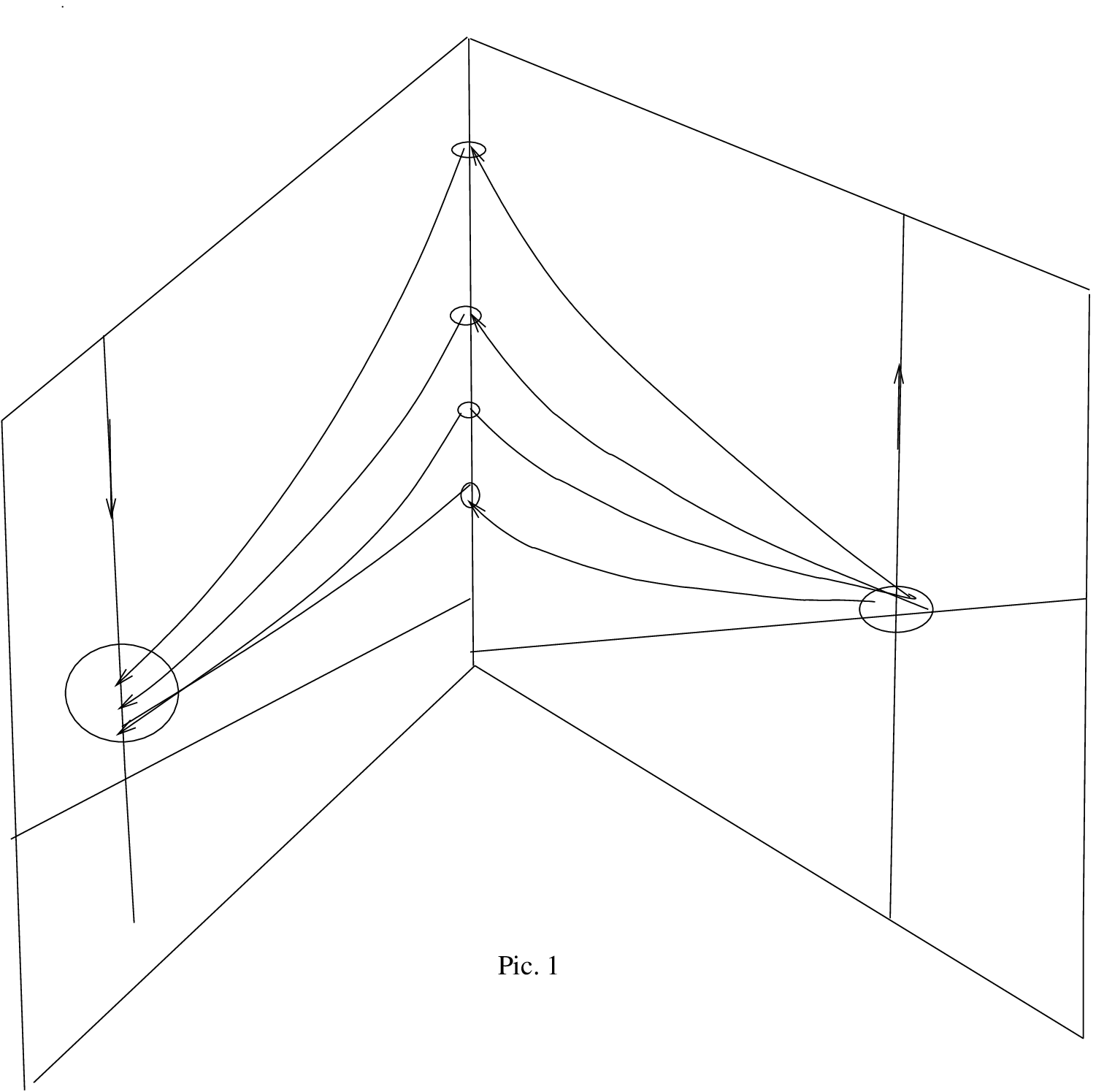} \vspace{15pt}
   \caption {Transversal pairs }
   \end{figure}
\noindent\pro{\it Proof}\endpro We can and will assume that the solvable radical of $G$ is unipotent. Indeed let $\Gamma_1 = [\Gamma, \Gamma]$ and $G_1 =[G,G].$ Let $R$ be the solvable radical of $G.$ It is well known that $[G,R] \subseteq U.$ Hence the solvable radical of $G_1$ is unipotent. Obviously $G_1$ is the Zariski closure of $l(\Gamma_1)$ and fulfills all requirements of the proposition. Thus if $\Gamma_1$ does not act properly discontinuously then the same is true for $\Gamma.$  \\
Let $\tilde{S}$ be a maximal reductive subgroup of a connected group $G$ whose solvable radical $R$ is unipotent. Then $\tilde{S} \cap R =\{1 \}$. Thus $S$ is a maximal reductive subgroup of $G$. Consequently, every regular element $g$ of $G$ is conjugate to an element of $S.$ Let $\sigma: G \rightarrow S$ be the projection.  The set of regular elements in $S$ is Zariski open. Since $\Gamma$ is Zariski dense in $G$ there exists an element $g \in \Gamma$ such that $\sigma (l(g))$ is a regular element of $S$. Let $x =l(g)$ and let $x =x_s x_u$ be the Jordan decomposition of $x.$ Thus $\sigma(x)$ is a regular element of $S.$ Therefore $\sigma (x_u)=1$ and consequently $x_u \in U.$ By the arguments above, $x_s$ is conjugate to an element of $S.$ Hence we can and will assume that $x_s \in S.$ As a result we have $\sigma (x) =x_s$ and $x_u \in U.$ From (1) and(2) follows that $l(g) =x_s \in S.$ Indeed, by (1), $x_u v -v \in V_0$ for every $v \in V.$ By direct calculations from $x_s x_u=x_u x_s$ and (2) we conclude that $x_u =1.$ Hence, $l(g)$ is a regular element in $S$. By (3) the element $l(g)$ can be transformed into a transversal pair inside $G$. The set $\Omega_{ab}(l(g))$ is clearly not empty. By Proposition 2.6, (A1), (A2) and (A3)  the element $l(g) \in l(\Gamma)$  can be transformed into a transversal pair by an element $l(t), t \in\Gamma,$ such that the elements $g$ and $tgt^{-1}$ do not commute.\\
Set $h = t g^{-1} t^{-1}.$
 Clearly $A^-(h) = l(t)A^+(g).$ Since $t_g \in V_0$ it  follows from (1) and (**)  that $t_h = -t_g.$ In particular the lines $L_g$ and $L_h$  are parallel. Set $v = t_g.$ \\
By definition 2.5  there exists a 
subspace $W \subseteq  A^+(g)$ such that  $ l(t)W \oplus D^+(g) =V.$
Put $\tilde {W} =l(t) W.$ Clearly, $\tilde {W} \subset A^-(h).$ Obviously, the intersection 
 $(L_h +\tilde{W})  \cap E^+_{g} =L$ is a one dimensional affine space. Moreover, since $L_g$ and $L_h$ are parallel,
 $L$ is parallel to each of them. Since $h$ and $g$ do not commute we conclude that $L_g \cap L_h = \varnothing.$ Otherwise $\Gamma$ does not act properly discontinuously. There exists a constant $c =c(g,h)$ such that the distance\\ $d(L_g,L)\leq c.$ 
\par Fix a point $p_1 \in L_g.$   There exists a point $p\in L$ such that the vector $\overrightarrow{p p_1}$  is in $D^+(g).$ Let $p_2$ be a point in $L_h$ such that $\overrightarrow{p p_2} \in \tilde{W}.$ Let $U_d(p_1)$ be the ball in $D^+(g)$ of radius $d$ with the center at $p_1$   and  let $U_a (p_2)$ be the ball in $L_h+ A^+(h)$ of radius $a$ with the center at $p_2.$ We can and will assume that $ U_d (p_1) \cap  U_a (p_2) =\varnothing.$ It is easy to see, that there exist $N \in \mathbb{Z}, N >0,$ such that for every point $x_n = p + n v$  we have $g^{-n} x_n \in U_d(p_1)$ and $h^{n} x_n \in U_a(p_2)$ for every $n > N$ (see Pic.1). Thus for every $n > N$ there exists a point $y_n \in U_d(p_1)$ such that $h^{n}g^n y_n \in U_a(p_2).$ Hence $h^{n}g^n \neq  1$ and $h^{n}g^n U_d(p_1) \cap U_a(p_2) \neq \varnothing$ for all $n >N, n \in \mathbb{N}.$ Thus
\begin{enumerate}
\item[(1)]  $h^{Nm_1}g^{Nm_1} \neq  h^{Nm_2}g^{Nm_2} $  for all $m_1 \neq m_2, m_1, m_2 \in \mathbb{N} .$
\item[(2)]  $h^{Nm}g^{Nm} U_d (p_1) \cap  U_a(p_2) \neq \varnothing$ for all $m \in \mathbb{N}.$
\end{enumerate}
 Therefore the group $\Gamma$ does not act properly discontinuously. \\
We will prove a slightly more general statement.\\ 
\noindent \pro {\textit{Proposition 2.10}} Let $G \subset GL(V) $ be the Zariski closure of the linear part of an affine group $\Gamma.$ Let $S$ be a semisimple part of $G$ and let $U$ be the unipotent radical of $G.$ Assume there exists a chain of  $l(G)$-invariant subspaces ${0} \subseteq V_0 \subset V_1 \subseteq V_2=V$ such that the following conditions  hold.
\begin{enumerate}
\item[(1)] the induced representations of $S$ on $V_2 /V_1$ and $V_0$ and the induced representation of $U$ on $V_1/V_0$ are trivial. 
\item[(2)] Let $i : S \rightarrow SL(V_1 /V_0)$ be the induced representation of $S$. Then for one (then for every) regular element $g$ of $S$ the element  $i(g)$ does not have one as an eigenvalue.
\item[(3)] Every regular element $s \in S$  can be transformed into a transversal pair inside $S$. 
\end{enumerate}
Then $\Gamma$ does not act properly discontinuously.
\endpro\\
\noindent\pro {\it Proof.}\endpro Let $\Gamma_1 =[\Gamma, \Gamma]$ and let $G_1$ be the Zariski closure of $\Gamma_1.$ 
From (1) follows that the solvable radical $R$ of $l(G_1)$ is a unipotent subgroup of $GL(V).$ Let
 $\Gamma_{m+1} =[\Gamma_1,\Gamma_m ]$ and let $G_m$ be the Zariski closure of $\Gamma_m,  m \geq 1, m \in \mathbb{Z}.$ It is well known that $G_{m+1} = [G_1, G_m].$\\
\noindent There exists an $N \in \mathbb{N}$ such that for all $ m \geq N, m \in \mathbb{Z},$ the restriction of $l(G_m)$ to $V_0$ and the induced action of $l(G_m)$ on $V_2/V_1$ are trivial. Since a semisimple part of $G$ is also a semisimple part of $G_m$ for all $m \in \mathbb{N}$ we conclude that $\Gamma_m$ fulfils all requirements of the proposition. Assume that  $m \in \mathbb{Z}, m> N.$ We will show that the group $\Gamma_m$ does not act properly discontinuously.\\
Indeed, since the induced representation $l(G_m)$ on $V_2/V_1$ is trivial it follows from 2.1 $(\ast)$,  that the affine subspace $\overline{V}_1= V_1 +0$ is $G_{m+1}$ invariant. Denote by $\overline{\Gamma}$ (resp. $\overline{G}$ ) the restriction of $\Gamma_{m+1}$ (resp. $G_{m+1}$) to
$\overline{V}.$ 
\par  If $V_0 ={0}$ then for every regular element $\g \in \Gamma_{m+1}$ there exists a fixed point $q_{\gamma}.$  Hence  $\Gamma_{m+1} $ does not act properly discontinuously. Since $\Gamma_{m+1} \leq \Gamma$ the group $\Gamma$ does not act properly discontinuously.   Assume that $V_0 \neq {0}.$  Obviously  $\overline{\Gamma}$  and   $\overline{G}$ fulfil the hypotheses of Proposition 2.9. Hence  $\overline{\Gamma}$  does not act properly discontinuously.  Hence by the same argument as above we conclude that $\Gamma$ does not act properly discontinuously. This proves the proposition.

\section{ Possible linear parts}\label{0}
\pro {\it 3.0. Notation and terminology.}\endpro Let $\G$ be an affine group acting properly discontinuously. Let $G $ be the Zariski closure of $\G$ and let $S$ be a semisimple part of $G$. Clearly, $S$ is a semisimple part of the connected component of the linear part $l(G)$ of $G$. The goal of this section is to give a complete list of all possible semisimple subgroups $S$ of  $GL(V)$, $V= \mathbb{R}^n,$ which might be a semisimple part of $l(G)$. 
The possible semisimple subgroups of $GL(V)$,  which occur in our list fulfil the following assumptions (A1) and (A2). 
\begin{enumerate}
\item[(A1)] $dim V \leq 6$. 
\item[(A2)]\textit{There is a simple normal subgroup } $S_1 \leq S$  \textit{with} $rank_{\mathbb{R}}(S_1) \geq 2.$
\end{enumerate} 
The justification for (2) comes later in the beginning of the next chapter 4.\\
By Proposition 2.8 every element of the connected component $l(G)^0$ of $l(G)$ has one as an eigenvalue. Therefore we add to our assumptions the following one. 
\begin{enumerate}
\item[(A3)]\textit{Every  element $g \in l(G)^0$ has one as an eigenvalue}.
\end{enumerate}
\par We call a linear algebraic group $G \leq G_n$  with semisimple part $S$
 an \textbf{A}--group if they fulfil the assumptions (A1), (A2) and (A3). 
\par The main steps to establish our list are the following. For a semisimple group $S$ satisfying the properties (A1) -(A3) we shall see that there are at most two non-trivial irreducible components $V_i , i \leq 2,$ of the representation of the complexification $\overline{S}$ of $S$ on the complexification  $\overline{V}$ of $V$ and that the image $\overline{S}_i$ of $\overline{S}$ in $GL(\overline{V_i})$ is a simple group for every non-trivial irreducible component $\overline{V}_i,$ see 3.4. Furthermore it does not happen that $\overline{V}$ contains two non-trivial irreducible components $\overline{V}_1$ and $\overline{V}_2$ such that $\overline{S}_1$ and $\overline{S}_2$ are isomorphic, see 3.2. It follows, that  if the real Lie group $S$ is simple then also $\overline{S}$ is simple, see 3.3. Note that there are several ways to satisfy (A3). If for every non-trivial irreducible component $W$ of the representation of $S$ on $V$ there is an element $s \in S$ with no eigenvalue one on $W$ then the space $V_0$ of $S$-fixed vectors in $V$ has positive dimension, in particular the sum of the dimensions of these spaces $W$ is at most 5.
\par We will assume from now on that $G$ is an \textbf{A}--group. If the dimension of every simple normal subgroup of a semisimple part $ S$ is $\leq 6$ then $(A2)$ does not hold. Thus there exists a simple normal subgroup of $S$ with dimension $ > 6.$ 
 Let us now recall a list  [PV, pp  260-261] of all possible complex representations   $\rho$  of a simple complex Lie group $S$ with  $\dim \rho \leq 6 \leq  \dim S.$   In the first column the symbols $SL_n,\, Sp_{2n},\, SO_n$ denote the corresponding simple Lie (algebraic) group in their simplest representation. The symbol $S^m H$ (resp. $\wedge^m H$) denotes the $m^{th}$ symmetric 
(resp. exterior) power of a linear group, and $S^m_0 H$ (resp. $\wedge^m_0 H$) is the highest (Cartan) irreducible component of this representation. 
\par \textbf{Table 1}\\
\medskip
\begin{tabular}{|l|l|l| }
\hline
$ S $ & $\dim \rho$ & $n$\\
\hline
$SL_n, n \geq 3$ & $n $ & $n =3,4,5$  \\
\hline
$SO_n, n\neq 4, n \geq 3 $& $n  $ &$n =3, 5, 6$\\
\hline
$Sp_{2n}$ &$2n$ &$2,3$\\
\hline
$AdSL_n$ & $n^2 -1 $ & $n = 2$\\
\hline
$S^2Sl_n$ & $n(n+1)/2 $& $n=2,3 $ \\
\hline
$\wedge^2 SL_n , n \geq 4$& $n(n-1)/2$&$n=4 $\\
\hline
$\wedge^2 SO_n , n \geq  3, n \neq 4 $& $n(n-1)/2$&$n=3$ \\
\hline
$\wedge_0^2 Sp_{2n} , n \geq 2$& $(n-1)(2n+1)$ & $n=2 $ \\
\hline
\end{tabular}\\
\medskip\\
Our next goal is to provide a list of all possible real simple linear groups $S$ which might  be a semisimple part of $G$ and which are possible as a factor of a semisimple part of an \textbf{A}--group.
 We will use the following notation. Let $\overline{V} =V\otimes_{\mathbb{R}}\mathbb{C}$ be the complexification  of $V$ and let $\overline {S}$ be the complex Lie group, such that $S$ is a real form of $\overline {S}.$
\par If the group $\overline{S}$ is simple then the group $S$ is a simple real Lie group.  Assume that the space  $\overline{V}$ is irreducible. Then $\overline{S}$ is a group listed in Table 1.\\ 
\pro {\it 3.1. Simple and irreducible.}\endpro Thus using [OV] we have the following list of all real simple groups $S$ which are a real form of a simple complex Lie group $\overline {S}$ listed in the Table 1 no matter if they are \textbf{A}--groups or not.
\pagebreak
\par\textbf{Table 2 }\\
\medskip

\begin{tabular}{|l|l|}
\hline
$ S $ & $\dim \rho$ \\
\hline
$SL_k(\mathbb{R}),  3 \leq k \leq n$ & $k$  \\
\hline
$SO(3,2) $& $5$ \\
\hline $SO(2,1) $& $3$\\ 
\hline $SO_3(\mathbb{R})$& $3$\\
\hline $Sp_{4}(\mathbb{R}) $& $6$ \\
\hline 
\end{tabular}\\
\medskip\\
\pro {\it 3.2. Simple and reducible.}\endpro
  Here we assume that the space  $\overline{V}$ is reducible and the complex group  $\overline {S}$ is simple. There exists an $\overline {S}$-invariant non-trivial subspace $W \leq  \overline{V}$ with non-trivial representation of $\overline {S}.$ Obviously $\dim W \geq 2.$ If $\dim W = 2$ then \textit{rank}$S \leq 1.$ If $\dim W =3$  the real form of $\overline{S}$ does not have one as an eigenvalue if \textit{rank}$S \geq 2.$ Thus there is no simple \textbf{A}--group such that $\overline{V}$ is $\overline{S}$ -reducible. \\
\pro {\it 3.3. Semisimple not simple}\endpro Let $S$ be a simple Lie group, such that $\overline {S}$ is not simple and $\overline{V}$ is irreducible.
There exists a complex structure on $S.$ Namely, there is a complex simple Lie group $\tilde{S}$ such that $S =\tilde{S}(\sigma(\mathbb{C}))$, where $\sigma : \mathbb{C} \longrightarrow M_2(\mathbb{R})$ is the natural embedding of the field $\mathbb{C},$ see [OV]. In this case  $\overline{S} $ considered as real Lie group is isomorphic to $ \tilde{S} \times \tilde{S}.$ The possible groups of this type are listed in the following table 3. \\
\par \textbf{Table 3 }\\
\medskip

\begin{tabular}{|l|l|}
\hline
$ S $ & $\dim \rho$ \\
\hline
$SL_k(\sigma(\mathbb{C})),  k =2, 3$ & $2k$  \\
\hline $SO_3(\sigma(\mathbb{C}))$& $6$\\
\hline 
\end{tabular}\\
\medskip\\
Non of the groups listed in Table 3 fulfils properties \textbf{A}. Indeed, since $\dim V \leq 6$ the groups $SL_2(\sigma(\mathbb{C}))$, $SO_3(\sigma(\mathbb{C}))$ have real rank 1 and $SL_3(\sigma(\mathbb{C}))$  does not have one as an eigenvalue. Note that non of the groups listed in Table 3 can be a normal subgroup of the semisimple part of $G$.\\ 
\pro {\it 3.4. General case. }\endpro The semisimple group $\overline{S}$ is the almost direct product of simple groups $\overline{S}= \prod_{1 \leq i \leq k} \overline{S}_i, k \geq 2 .$  Let $W_0 =\{v \in \overline{V}:   s v =v\,  \forall s \in \overline{S} \,\}.$ There exists a unique $\overline{S}$--invariant subspace  $\overline{W}$ of  the space $\overline{V} $ such that $\overline{V} $  is the direct sum of $W_0$  and  $\overline{W}.$  If the restriction  $\overline{S}\vert _{\overline{W}}$ is an irreducible representation of  $\overline{S}$, then it is the tensor product of $\overline{S}_i $- irreducible representations  for all $i =1, \dots, k.$ Thus if  $\dim \overline{V} \leq 6$ it follows immediately that this is impossible for an \textbf{A}--group. Therefore  $\overline{W}$  is the direct sum of $\overline{S}$--invariant non-trivial irreducible subspaces  $W_i,  i =1, \dots, k$ such the restriction $ \overline{S_j} \vert _{W_i}$ is trivial for every $i \neq j, i,j =1, \dots, k.$   As we know, every element of $G^0$ has one as an eigenvalue. Thus it follows from $(A3)$ that  if the subspace $W_0$ is trivial, there exists an $i_0, 1 \leq i_o \leq k,$ such that every element $s \in \overline{S}_{i_0}$ has one as an eigenvalue.  Since for every $i =1, \dots, k$ the group $\overline{S}_i$ is an irreducible subgroup of  $GL(W_i )$,  we can and will again use Table 1 and Table 2. This will lead us to a complete list of all possible cases under the assumption \textbf{A}. \\
\noindent\pro{\it 3.5. Linear parts and decomposition.}\endpro Let us summarize all  we did in 3.1, 3.2, 3.3 and list all  cases we have to consider.  Let $V_0$ be the maximal subspace in $V =\mathbb{R}^n$ such that $S$ acts  trivially on $V_0$. Let $V_1$ be the unique $S$--invariant subspace such that $\mathbb{R}^n = V_0 \oplus V_1.$
Let $\pi_S : G  \longrightarrow S$ be the projection. We will use these notations throughout the rest of the paper. Recall that $G$ is an \textbf{A}--group.
\bigskip \\
\pro {Case 1} \endpro Assume that for every regular element $s  \in S$ the restriction $s|_{V_1} $ does not have $1$ as an eigenvalue. Thus $V_0 \neq {0}$. Consider the inclusion\\ $i_s: S \longrightarrow GL(V_1)$ as a representation of the semisimple Lie group $S.$ Assume first $S$ is a simple group.It follows from 3.3 that the complexification of $i_s(S)$ is a simple irreducible group. Thus it follows from the Table 2 that  all possible semisimple parts which have property (A2)  are:
\begin{enumerate} \item[(1)] $S= SL_l(\mathbb {R}),  V_1 =\mathbb{R}^l,  l < n,  3 \leq l \leq 5,  4 \leq n \leq 6.$
\item[(2)] $S= Sp_4(\mathbb{R}), V_1 =\mathbb{R}^4, n=5, 6.$
\end{enumerate}
Suppose that the group $S$ is semisimple, but not simple. As we show  in this case $i_s(S)$ is the direct product of two simple groups such that their complexifications are simple complex groups. It follows from Table 2  that  all possible semisimple parts in this case which have property (A2) are:
\begin{enumerate}
\item[(3)] $S= SL_2(\mathbb{R}) \times SL_3(\mathbb{R}),  V_1 =\mathbb{R}^5, n=6.$
\end{enumerate}
\pro {Case 2.}\endpro  Assume that for every regular element $s  \in S$ the restriction $s|_{V_1} $ has $1$ as an eigenvalue. Suppose that $S$ is a simple group. It follows from 3.2 and 3.3  and Table 2  that the group $\overline{S}$ is simple. Therefore $S =SO(3,2)$ and $\dim V =5, 6.$  If $S$ is a semisimple but not a simple group, we show above (see 3.4) that $S$ is the almost direct product of two simple group $S_1$ and $S_2$ such that their complexifications $\overline{S}_1$ and $\overline{S}_2$ are simple complex groups. Since $S$ is an \textbf{A}--group  it follows from Table 2 and (A2) and (A3)  that  $\dim V =6$ and $S= SL_3(\mathbb{R}) \times SO(2,1),$ or $ S= SL_3(\mathbb{R}) \times SO(3).$ Therefore we conclude that in this case
\begin{enumerate}
\item[(1)] $S =SO(3,2), \dim V_1= 5, n=5,  6.$
\item[(2)] $S =SO(3)\times SL_3(\mathbb{R}),   V_0 =0,   n = 6.$
\item[(3)] $S =SO(2,1)\times SL_3(\mathbb{R}),  V_0=0,  n=6.$
\end{enumerate}
Case 1 and 2 give us the complete list of all possible semisimple parts of the \textbf{A}--group $G$.\\
\pro{Remark 3.6.  }\endpro We see from the lists above that if $\dim V \leq 5$ a semisimple part of an \textbf{A}--group is a simple group.
 \section{The  Auslander conjecture in dimensions 4 and 5}\label{0}
\par \pro {\it 4.0.}\endpro\,Let us first explain the plan of the proof of the Main Theorem.
Let $\G$ be a crystallographic group and let $G$ be the Zariski closure of $\Gamma$. Then we have the Levi decomposition $G =S R$ where  $R$ is the solvable radical and $S$ is a semisimple part of $G$. Let $S = \prod_{1 \leq i \leq k} S_k $ be the decomposition of the semisimple part into an  almost direct product of simple groups. It is well known that if $rank _{\mathbb{R}}(S_i) \leq 1$ for all
$1 \leq i \leq k $ then $\G$ is not crystallographic [S2], [To2]. Therefore from now on unless otherwise noted  we will assume that in case $S \neq {1}$ we have $$\max _{ 1 \leq i \leq k} rank_{\mathbb{R}}(S_i) \geq 2.$$ Hence
$l(G)$ has property $(A2).$ This is the justification for our assumption $(A 2)$. 
By Proposition 2.8  we conclude that  $l(G)$ has property $A(3).$ Hence $l(G)$ is an \textbf{A}--group for $\dim V \leq 6.$ Therefore if the semisimple part $S$ is non-trivial it is one of the groups listed in Case 1 and Case 2. Our strategy is to show case by case that none of the semisimple groups listed in Case 1 and Case 2 is a semisimple part of the Zariski closure of a crystallographic group $\Gamma$.  Thus, $S={1}$ and $\Gamma$ is virtually solvable. \\ 
\pro {\it 4.1.}\endpro  In this section we will prove the Auslander conjecture in dimensions  4 and 5.\\
 Let $\Gamma$ be a discrete subgroup of an affine group and let $G$ be the Zariski closure of $\Gamma.$  Assume that the linear part $l(G)$ of $G$ is a connected  \textbf{A}--group and the semisimple part $S$ of $G$ is a simple group . 
Let us now recall the list of all possible cases for $n=\dim V \leq 6.$ 
It follows from 3.5 that all possible cases are $$  S= SL_l(\mathbb {R}),  V_1 =\mathbb{R}^l, 3 \leq l \leq 5, l < n,  4 \leq n \leq 6 \eqno (s_1)$$
$$ S= Sp_4(\mathbb{R}), V_1 =\mathbb{R}^4, SO(3,2), n=5,  6 \eqno (s_2)$$
We will deal with the case $S =SO(3,2), \dim V =6$, in the next chapter. Therefore in the next proposition we will assume that if $S=SO(3,2)$ then $\dim V =5.$\\
\noindent \pro{\it Proposition 4.2.} Let $\Gamma$ be a discrete subgroup of an affine group and let $G$ be the Zariski closure of $\Gamma.$ Assume that the simple part $S$ of $G$ is as  in
$(s_1)$ or $(s_2).$ Then the group $\Gamma$ does not act properly discontinuously.\endpro.\\
\pro {\it Proof.} \endpro 
\pro{ Case 1}\endpro$S = SO(3,2), \dim V =5.$ By Theorem B [AMS3] the group $\G$ does not act properly discontinuously.\\
\pro{ Case 2} \endpro$S \neq  SO(3, 2)$. 
Obviously $S$ fulfils all requirements of Proposition 2.10.
Thus $\Gamma$ does not act properly discontinuously. This proves the proposition.\\
 \pro {\it Proposition 4.3.} Let $\G  \subseteq \Aff(\mathbb{R}^n), n \leq 5$  be a crystallographic group. Then $\G$ is virtually solvable. \endpro\\
\pro{\it Proof}\endpro  Let $G$ be the Zariski closure of $\Gamma.$ Since the connected component $G^0$ is a finite index subgroup of $G$ we conclude that $\Gamma \cap G^0$ is a finite index subgroup of $\Gamma.$ Clearly $\Gamma \cap G^0$ is a crystallographic group. Thus we shall and will assume that the group $G$ is connected. As we explained above in 4.0, the group $l(G)$ is an \textbf{A}--group. Assume that $l(G)$ has a non-trivial semisimple part $S.$ It follows from 3.6 that for  $n \leq 5$  the group $S$ is a simple group listed in $(s_1)$ or $(s_2).$ Thus by Proposition 4.1 $\Gamma$ does not act properly discontinuously. Therefore, the semisimple part $S$ is trivial. Hence the group $\Gamma$ is virtually solvable. 
\section{The Auslander conjecture in dimension $6$. The cohomological argument.}\label{0}
\par We start use the same notations as in Chapter 4. The goal of this chapter is to show that if $\Gamma$ is a crystallographic group and $\dim V =6$ then the semisimple part of the Zariksi closure of $\Gamma$ can not be one of the groups listed in Case 1  and Case 2 (1), (2). The groups of Case 2 (3) will be dealt with in the next section 6.
We will with the following \\
\pro {\it Proposition 5.1.} Let $\Gamma$ be an affine group and let $G$ be the Zariski closure of $\Gamma.$  Assume that the semisimple part $S$ of $G$ is as in the Case 1 (1), (2). Then the group  $\G$ does not act properly discontinuously \endpro\\
\pro{\it Proof.}\endpro The proof follows immediately from Proposition 2.10.
\endpro \\
\pro {\it Proposition 5.2} Let $\Gamma$ be an affine group and let $G$ be the Zariski closure of $\Gamma.$  Assume that the semisimple part $S$ of $G$ is as in Case 1 (3)  $S=SL_2(\mathbb{R}) \times SL_3(\mathbb{R}).$  Then the group $\Gamma$ does not act properly discontinuously.\endpro \\
\pro {\it Proof.} \endpro We have a chain $0 \subseteq W_0 \subset W_1 \subseteq W_2 =V$ of
$l(G)$-invariant subspaces. There are three possible cases
\begin{enumerate}
\item[(i)] $\dim W_0 =1,$
\item[(ii)] $\dim W_1/W_0 =1,$
\item[(iii)] $\dim W_2/W_1 =1.$
\end{enumerate}
\pro{\it Cases (i) and (iii).}\endpro It follows from (A3) that in case (i) we have\\
$l(G) | _{W_0} =1$ and $\dim W_2/W_1 =0.$ In case (iii) the induced representation $l(G) \rightarrow GL(W_2/W_1)$ is trivial and $\dim W_0 =0.$ Hence $\Gamma$ does not act properly discontinuously by  Proposition 2.10 .\\
\pro{\it Cases (ii).}\endpro The induced representation $l(G) \rightarrow GL(W_1/W_0)$ is trivial as follows again from (A3). Roughly speaking the space of $S$-fixed vectors is "in between". Set $U_0 =W_0.$ There exist $S$--invariant spaces $U_1$ and $U_2$  such that $ W_1 = U_0 \oplus U_1 $,  and $V = U_0 \oplus U_1 \oplus U_2$,\\
We will prove the statement of the proposition assuming that $ S|_ {U_0} =SL_3(\mathbb{R}),$ $S|_{U_1} =I $ and $S|_{U_2} =SL_2(\mathbb{R}).$  The proof in  case $ S|_{W_0} =SL_2(\mathbb{R}),$ $S|_{U_1} =I $ and $S|_{U_2} =SL_3(\mathbb{R})$ 
is a verbatim repetition.\\
There exists a $g\in \Gamma$ such that $l(g)$ is an $\mathbb{R}$--regular element in $l(G)$ ([AMS1], [P]). We can and will assume that $l(g) \in S.$  Let $g_0 = l(g)|_{U_0} \in SL_3(\mathbb{R})$,  $g_1 = l(g)|_{U_1} =1$ and $g_2 = l(g)|_{U_2} \in SL_2(\mathbb{R}).$
We can assume that $ \dim A^-(g_0) < \dim A^+(g_0).$ Thus $\dim A^+(g_0)=2.$  Note that  $ \dim A^+(g_2)= 1$ and $A^0(g) = U_1.$ Let $U$ be a one dimensional $l(g)$--invariant subspace of $A^+(g_0).$ Then there exists $t \in S$ such that $l(t)U \notin A^+(g_0) \cup A^-(g_0) $ and $l(t)A^+(g_2) \notin A^+(g_2) \cup A^-(g_2)$ and  $l(t)U \oplus A^+(g_0) =U_0$ and  $l(t)A^+(g_2) \oplus A^+(g_2)  = U_2.$  Set  $A(t)= l(t)U +l(t)A^+(g_2).$ Then $A(t)  \oplus D^+(g)= V$ since $D^+(g)= A^+(g_0) +U_1 + A^+(g_2) $  Let $ \sigma : G \rightarrow S$ be the projection. Clearly $\sigma (\Gamma)$ is Zariski dense in $S.$ Therefore we can and will assume the $t \in \Gamma.$ Put $h = t g^{-1}t^{-1} \in \Gamma.$  Clearly $A(t) \subseteq A^-(h).$ Remark, that $0 \neq u \in U$ is an eigenvector of $h$ but not an eigenvector of  $g.$ Therefore $h^n \neq g^m$ for all $n,m \in \mathbb{Z}, n \neq 0, m\neq 0.$ 
 Let $A = l(t)U_1 + A(t),  A \subseteq  D^-(h) $ and let  $D =A +L_h.$ Clearly, $D$, is an $h$--invariant affine space in  $E^-_h$ and $\dim D \cap E^+_g =1.$ Let $L= D \cap E^+_g.$  We have the projections $\pi_1 : \mathbb{R}^6  \longrightarrow L_g $ of an affine space $\mathbb{R}^6$ onto $L_g$ along $A^+(g)+A^-(g) $ and $\pi_2 :\mathbb{R}^6  \longrightarrow L_h $ along  $A^+(h)+A^-(h) .$ The restriction $\overline{\pi}_i= \pi_i |_L, i =1,2$ is an affine isomorphism.  Set  $\theta= \overline{\pi}_2^{-1} \circ\overline{\pi}_1$. Then  $\theta : L_g \longrightarrow L_h$ is an affine isomorphism.   Since $g_1 =1$ we conclude $l(t) t_g -t_g \in U_0$. Combining this with  (**) we obtain $\theta (t_g) = -t_h.$
\par Let $p_1 \in L_g$ and let $p_2 \in L_h.$ There exists a point $p \in L$ such that  the vector $\overrightarrow{pp_1} \in A^+(g)$ and $\overrightarrow{pp_2} \in A^-(h).$ Consider a ball $U_1(p_1)$of radius $1$ and the center at $p_1$ and a ball $U_1(p_2)$of radius $1$ and the center at $p_2.$
there exists a point $q \in L $ such that the vector $\pi_1(\overrightarrow{pq}) =t_g .$ Set $x_k = p + kv, k \in \mathbb{N}, k > 0.$  Then there exists a positive $N,  N \in \mathbb{Z}$ such that for
$m > N$ we have $g^{-m}  x_m \in U_1(p_1)$ and  $h^{m}  x_m \in U_1(p_2).$ As in the proof of Proposition 2.9 we conclude that for $m > N$ we have $h^m g^{m}  U_1(p_1) \cap U_1(p_2) \neq \varnothing.$ Since $h^n \neq g^m$ for all $n,m \in \mathbb{Z}, n \neq 0, m\neq 0$   the group $\Gamma$ does not act properly discontinuously. \bigskip\\
\noindent \pro {Proposition 5.3. }  Let $\Gamma$ be an affine group and let $G$ be the Zariski closure of $\Gamma.$  Assume that $G$ is connected and the semisimple part $S$ of $G$ is as in Case 2 (1), (2). Then the group $\G$ is not a crystallographic group. \endpro \\ 
\pro {\it Proof} \endpro. Let us first explain the main idea of the proof. Since  the subgroup $\G \subseteq G_n$ is a crystallographic group, the virtual cohomological dimension $vcd(\G)$ of $\G$  is $\dim \mathbb{R}^n =n .$ Hence $vcd(\G)= 6.$ As a first step we will show that $vcd(\G) \leq \dim(S/K)$, where $S/K$ is the symmetric space of $S$. Then we compare $\dim S/K$ and $vcd(\G)$ in the cases $S =SO(3)\times SL_3(\mathbb{R})$, $S = SO(3,2)$ and come to the conclusion that $\dim S/K \geq vcd(\G).$ This will lead to a contradiction. We actually show that the projection $p:G \rightarrow S$ restricts to an isomorphism of $\Gamma$ onto a discrete subgroup of $S$. In case $S = SO(3) \times SL_3(\mathbb{R})$ the dimension of $S/K$ is 5 and so $vcd(\Gamma) \leq dim(S/K)$ is impossible. In case $S = SO(3,2)$ the dimension of $S/K$ is 6 and so $p(\Gamma)$ would be a cocompact lattice in $S$ and we will get a contradiction using the Margulis rigidity theorem.
\par Let us first show that $vcd(\G) \leq \dim(S/K).$ Let $R$ be the solvable radical and  $U$ be the unipotent radical of $G.$ Recall that $G$ acts trivially on the factor-group $R /U$. Thus it is easy to see that  in Case 2 (2) we have $R =U.$ Let $\G_r = R \cap \G$ and let $R_1$ be the Zariski closure of 
$\G_r.$ Then the group $R_1$ is a normal subgroup in $G$ since $\G_r$ is a normal  subgroup in $\G.$  Let $T_1$ be a maximal torus and let $U_1$ be the unipotent radical of $R_1$. Since $\tilde{S} =S T_1$ is a reductive subgroup of $G$ there exists a point $q_0$ such that $\tilde{S} q_0 =q_0$ [see (2.1)]. 
Set $W =R_1 q_0$. For every $g \in R_1$ there are unique elements $ t \in T_1$ and $u\in U_1$ such that $g =t u.$ Define a map $\pi : R_1 \rightarrow U$ by
$\pi (g) =u.$ Then $\pi (\G_r)$ contains a uniform lattice of $U_1$ [S2, Proposition 2]. Since $T_1 q_0 =q_0$ we conclude that $\G_r \setminus W$ is compact. 
\par Obviously $s W = W$ for $s \in S$, since $sq_0=q_0$ and $R_1$ is a normal subgroup of $G.$  Let $\rho: S \rightarrow GL(T_{q_0})$ be the representation of $S$ on the tangent space $T_{q_0}$ of $W$ at $q_0.$ It is clear that the only possible numbers for $\dim (T_{q_0})$  are $\{0, 3,  6\}$ if $S =SO(3)\times SL_3(\mathbb{R})$ and $\{0, 1, 5, 6 \}$ if $S= SO(3,2).$ Let us show that in each case $\dim(T_{q_0})=0.$ Assume that $\dim (T_{q_0})=6.$
Then $W =R_1q_0=\mathbb{R}^6.$ As we show above $\G_r \setminus W$ is compact. Thus $\G_r$ is a crystallographic group. On the other hand $\G_r$ is a subgroup of a crystallographic group $\G$ which acts on the same affine space. Thus the index $|\G / \G_r|$ is finite. On the other hand the index $|\G / \G_r|$ is infinite. Otherwise the Zariski closure of the solvable group $\Gamma_r$ would contain the connected component of the Zariski closure of  $\Gamma$ which is impossible. We thus have shown that $\dim(T_{q_0}) < 6.$ 
\par We will treat the two cases $S=SO(3) \times SL_3(\mathbb{R})$ and $SO(3,2)$ separately.
\par Let $S =SO(3) \times SL_3(\mathbb{R})$ and $\dim (T_{q_0})=3.$ Then $G$ is a subgroup of the following group $\widetilde{G} = \{ X : X \in GL_7(\mathbb{R}) \},$  where
$$X=
\left(
\begin{array}{ccc}
A&B& v_1\\
0&C&v_2\\
0&0&1\\
\end{array}
\right)\,
 \eqno (1) $$ 
where\,  $A \in SO(3), C \in SL_3(\mathbb{R}), v_1, v_2 \in \mathbb{R}^3, $\\
or 
$$X=
\left(
\begin{array}{ccc}
A&B& v_1\\
0&C&v_2\\
0&0&1\\
\end{array}
\right)
\eqno (2) $$
 where \,$ A \in SL_3(\mathbb{R}) , C \in SO(3) , B \in M_3(\mathbb{R}), v_1, v_2 \in \mathbb{R}^3.$\\
We claim that $\dim T_{q_0} =0.$ We will prove this for  (1). The proof for (2) will go along the same lines. 
\par Since the semisimple part of a group has to commute with one maximal reductive subgroup of its solvable radical the solvable radical of $G$ is unipotent. Therefore for a $X \in R$ we have
$$X=
\left(
\begin{array}{ccc}
I_3&B(X)& v_1(X)\\
0&I_3&v_2(X)\\
0&0&1\\
\end{array}
\right)
.$$
Assume that there is an element $X$ of the unipoten $U$ radical of $G$ such that $B(X) \neq 0.$ Since $l(U)$ is a normal subgroup of $l(G)$ by direct calculations we show that for every 
$B \in M_3(\mathbb{R})$ there exists $X \in U$ such that
$$l(X)=
\left(
\begin{array}{cc}
I_3&B\\
0&I_3\\
\end{array}
\right)\,
\eqno (3). $$
Otherwise $B(X) =0$ for every element $X \in U.$
\par The unipotent group $R_1$ is a normal connected subgroup of $G$ and $R_1 \leq U.$ There are three connected proper nontrivial normal unipotent subgroups of $G,$ namely, $R_1 = \left\lbrace X \in U, v_2(X) =0 \right\rbrace ,$  $R_1 = \left\lbrace X \in U, B(X) =0,  v_2(X) =0 \right\rbrace $ and
 $R_1 = \left\lbrace X \in U, B(X) =0,  v_1(X) =0 \right\rbrace. $
We conclude that in these cases $W =R_1q_0$ is an affine $G$--invariant subspace.  Thus we have a nontrivial $G$--invariant affine space $W$ where  $\Gamma$ and $\Gamma_r$ act as crystallographic groups. By the same argument we used in case $\dim T_{q_0} =6$ we conclude that the subgroup $R_1 $ is trivial.
By Auslander's theorem [R],  $\pi_S (\G)$ is a discrete subgroup of $S.$ Since the intersection  $\G \cap R$ is trivial,  $\pi_S (\G)$ and $\G$ are isomorphic. Hence  $vcd(\G) = vcd( \pi_S(\G)) \leq \dim S/K $, where $K$ is a maximal compact subgroup in $S.$ Thus  $vcd(\G) \leq 5.$ On the other hand, $vcd(\G) =6,$ a contradiction.
\par Let us now show that Case 2  (1) is also impossible. We will use the notation introduced in 3.5. Recall that $V =V_0 \oplus V_1$ where the restriction $S |_{V_0} $ gives a trivial representation and the restriction
$S |_{V_1} =SO(3,2).$ Assume that $V_1$ is $l(G)$-invariant. Then it follows from the linear representation(*) in 2.1, that the affine space $V_1 + q_0$ is $\Gamma_1$ -invariant, where
 $\Gamma_1 =[\Gamma, \Gamma].$ Obviously  $\dim V_1 =5$ and $l(\Gamma_1) | _{V_1} \leq SO(3,2).$  It follows from Proposition 4.2, case 1, that $\Gamma_1$  does not act properly discontinuously on $V_1+q_0.$ Therefore $\Gamma$ is not a crystallographic group. Thus we can and will assume that $V_0$ is $l(G)$ invariant. 
We will prove first that  $\dim W = 0$. Recall that $W =R_1 q_0$ and $l(G)q_0=q_0.$ We have the following matrix representation of $G$. Let $X \in G$ then 
$$X=
\left(
\begin{array}{ccc}
\lambda(X)&w(X)& a(X)\\
0& A(X)& v(X)\\
0&0&1\\
\end{array}
\right)
,$$
where $A(X) \in SO(3,2), w(X), v(X) \in \mathbb{R}^5, \lambda(X),  a(X) \in \mathbb{R}.$  
As we concluded above, there are three possible cases for $\dim W $, namely, $\dim W = 0, 1, 5.$ Our goal is to show $\dim W \neq  1, 5.$
\par Assume that $\dim W =1.$ The representation $\rho$ of $ S$ on $T_{q_0}$ is trivial. Clearly,  $S=SO(3,2)$ is an irreducible subgroup of $GL( \mathbb{R}^5).$
Therefore we conclude that if $X$ is an element in the normal subgroup $R_1$  of $G,$ then $v(X) =0.$ 
Thus for every $X \in R_1$ we have
$$X=
\left(
\begin{array}{ccc}
1&w& a\\
0&I_5& 0\\
0&0&1\\
\end{array}
\right)
,$$
where $w \in \mathbb{R}^5, a \in \mathbb{R}.$  Hence $W$ is an affine $\G$-- invariant subspace in $\mathbb{R}^6.$ Therefore we have a natural homomorphism $ \theta : \G \longrightarrow$ \text{Aff}$(\mathbb {R}^6/W).$  
 By [S2, Lemma 4], $\G/ \G_r = \theta(\G)$  is a crystallographic subgroup in \text{Aff}$(\mathbb R^6/W ).$ Obviously the semisimple part of the Zariski closure of $\theta (\G)$ is $SO(3,2)$ and $\mathbb R^6/W =\mathbb{R}^5.$ By [AMS3, Theorem A] this is  impossible. \\
  \par Assume that $\dim W =5.$ Again consider the space of orbits $\widehat{R} = \{ g W, g \in G\}$. Recall that the unipotent radical $U$ acts transitively  on $\widehat{R}$ [S2].  It is clear that $\widehat{R} $ is a one dimensional manifold. As in [S2] we have a representation $\rho$ of $l(G) $ on  the tangent space $T_W$  of $\widehat{R}$ at $W$. We show in [S2,  Theorem A] that one is an eigenvalue of  $ \rho (g)$ for every element  $g \in l(G).$ Hence the representation  $\rho$ is trivial. Note that this implies that $\lambda(X) =1$ for every $X \in R_1.$ Thus $R_1$  is a unipotent group.  Since  $\dim W =5$ there exists an $X \in R_1$ such that $v(X) \neq 0.$ On the other hand looking at the representation of $SO(3,2)$ on $R_1$ we conclude that if there exists $X \in R_1$ such that $v(X) \neq 0$ then for every $v \in \mathbb{R}^5$ there exists $X \in R_1$ such that $v(X) =v.$ 
Therefore for every $g \in 
\Gamma$ there exists $r \in R_1$ such that  for $\widehat{g} =g r^{-1}$ we have $v (\widehat{g}) =0.$ Obviously
 $$\widehat{g}=
\left(
\begin{array}{ccc}
1&w(\widehat{g})& a(\widehat{g}) \\
0&S(\widehat{g})& 0\\
0&0&1\\
\end{array}
\right)
  \eqno (4) .$$ 
 Let  $g_1, g_2 \in \Gamma$ be two elements on $\Gamma.$ There exist $ r_1, r_2 $ such that for
$\widehat{g}_i  r_i =g_i, i=1,2$ we have $v (\widehat{g}_i) =0.$
Clearly, we have $[\widehat{g_1}, \widehat{g_2}] \in l(G)$ and $[g_1, g_2] = [\widehat{g_1}, \widehat{g_2}]r_0 $ where $r_0 \in R_1$. Therefore $[g_1, g_2] W =  [\widehat{g_1}, \widehat{g_2}]r_0 R_1q_0 =
R_1 q_0=W. $ 
Set $\G/ \G_r = \theta(\G).$ By [S2, Lemma 4], $ \G/ \G_r =\theta(\G)$ acts as a crystallographic group on $\widehat{R}.$ Therefore $Stab_{\theta(\Gamma)}(W)$ is a finite set.  From $[\G, \G] W  =W$ and $\Gamma_r W =W$ follows that $[\theta(\G), \theta(\G)] \leq Stab_{\theta(\G)} (W).$ So  the group $[\theta(\G), \theta(\G)] $ is finite. Consequently $\theta(\G)$  is a virtually abelian group. Therefore $\Gamma$ is a virtually solvable group. This is a contradiction.
Thus we conclude that $W =0.$ Hence  $R_1 =\{e\}$ and the restriction of the homomorphism $\pi_S : G \longrightarrow  S =G /R $  onto $\G$ is an isomorphism. By Auslander's theorem [R], the projection $\pi_S(\G)$ is a discrete subgroup in $S$ and  $vcd(\pi_S(\G)) =vcd(\G)= 6.$ On the other hand $vcd(\pi_S(\G)) \leq \dim S/K$, where $K$ is a maximal subgroup in $S$. Obviously, $\dim S/K = 6.$ Hence $vcd (\pi_S(\G)) =\dim S/K.$ Therefore $\pi_S(\G)$ is a co-compact lattice in
$S.$ We can apply the Margulis rigidity theorem, since $rank_{\mathbb{R}} (S)=2$ and conclude that  there exists a $g \in \G$ such that $\tilde{\Gamma}= g \G g^{-1} \cap S $ is a subgroup of finite index in $\Gamma$. Since $\tilde{\Gamma} \leqslant S$ we have $\tilde{\Gamma} p_0 =p_0.$ Thus $\G$ does not act properly discontinuously.\\
\pro{\it Remark 2.} \endpro Under the assumption of Proposition 5.3 using the dynamical ideas and results from [AMS4] we can prove more.  Namely,  the group $\G$ does not act {\it properly discontinuously}. Hence in particular $\G$ is not crystallographic.  
\section{The Auslander conjecture in dimension $6$.  Dynamical arguments}\label{0}
\noindent\pro {\it 6.1. Orientation.}\endpro  The dynamical approach we have used in [AMS3] and will use here is based on the so called \textbf{Margulis sign} of an affine transformation.  The case $S =SO(2,1) \times SL_3(\mathbb{R})$ needs other tools, namely a new version of the Margulis sign. We will need to introduce it for the natural representation of $S$ which goes roughly saying by ignoring the $SL_3(\mathbb{R})$--factor. We then have a lemma similar to the cases of $SO(k+1,k)$, namely lemma 6.7, which says that if a group acts properly discontinuously, then opposite signs are impossible.
\par Now we will recall the important definition of the sign of an affine transformation. This definition was  first introduced by G. Margulis [M] for $n=3.$ Then it was generalized in [AMS3] for the  case in which the signature of the quadratic form is $(k+1, k)$ and finally for an arbitrary quadratic form in [AMS4]. Our presentation will follow along the lines of [AMS4]. Let $B$  be a quadratic form
of signature $(p,q),~p\geq q,~p+q= n$. Let $v$ be a vector in ${\mathbb{R}}^n,$ $v=
x_1 v_1+ \dots + x_p v_p +y_1 w_1 + \dots + y_q w_q$, where $v_1,
v_2, \dots, v_p,$\,  $w_1, w_2, \dots, w_q$ is a basis of
${\mathbb{R}}^n$. We can and will assume that
$$ B(v,v) = x^2_1 +\dots + x^2_p - y^2_1 - \dots - y^2_q. $$
Consider the set $\Phi$ of all maximal $B$\!--isotropic subspaces.
Let $X$ be the subspace spanned by $\{v_1, v_2, \dots, v_p\}$ and
$Y$ be the subspace spanned by $\{w_1, w_2, \dots, w_q\}$. It is
clear that ${\mathbb{R}}^n = X \oplus Y$. Define the cone
$$Con_{B}= \{v \in \mathbb{R}^n : B(v, v) < 0 \} .$$
Clearly $Y \backslash \{0\} \subset Con_{B}.$
We have the two projections
$$
\pi_X : \br^{n}\la X\ \rm { and } \,  \pi_Y: \br^{n}\la Y
$$
along $Y$ and $X$, respectively.  The restriction of $\pi_Y$ to
$W\in\Phi$ is a linear isomorphism $W\la Y$. Hence if we fix an
orientation on $Y$, then we have also fixed an orientation on each
$W\in\Phi$. For $W\in\Phi$, let us denote the $B$\!--orthogonal
subspace of $W$ by $W^{\perp} = \{z\in\br^{n}\ ;\ B(z,W) = 0\}$.
We have $W\subset W^{\perp}$ since $W$ is $B$\!--isotropic. We
also have
$$
\dim W^{\perp} = \dim W +(p-q)= p.
$$
The restriction of $\pi_X$ to $W^{\perp}$ is a linear isomorphism
$W^{\perp}\la X$. Hence if we fix an orientation on $X$, then we have
also fixed an orientation on $W^{\perp}$ for each $W\in\Phi$. Thus
we have orientations on both $W$ and $W^{\perp}$ and we have
naturally induced an orientation on any subspace $\widehat{W}$, such that
$W^{\perp}=W  \oplus \widehat{W} $. If $V_1\in\Phi$ and $V_2\in\Phi$ are
transversal, then $V_0 = V_1^{\perp}\cap V_2^{\perp}$ is a subspace
that is transversal to both $V_1$ and $V_2$; therefore $V_0 \oplus
V_1= V_1^{\perp}$ and $V_0 \oplus V_2=V_2^{\perp}$. So there are two
orientations $\o_1$ and $\o_2$ on $V_0$, where $\o_i$ is defined if
we consider $V_0$ as a subspace in $V_i^\perp$. We have [see AMS3, Lemma 2.1] \bigskip\\
\pro{\it Lemma 6.2.} The orientations defined above on $V_0$ are the
same if $q$ is even and they are opposite if $q$ is odd, i.e. $\omega_1 =(-1)^q \omega_2.$\endpro\\ Let us explain this in the special case when $p= k+1,  q=k .$\\
\noindent \pro{\it Example 6.3.} \endpro Let $V_1$ and $V_2$ be the maximal isotropic subspaces spanned by
the vectors $\{ w_1+v_1, \dots, w_k+v_k\}$  and $\{w_1-v_1, \dots, w_k-v_k\}$ respectively. Since for every $i=1, \dots, k$  we have $\pi_Y(w_i \pm v_i)) = w_i, i=1, \dots, k,$ we
conclude that  $w_1+v_1, \dots, w_k+v_k$ (resp. $w_1-v_1, \dots, w_k-v_k$) is a positively
oriented basis of $V_1$ (resp. $V_2$)
Then $V^\bot_1 \cap V^\bot_2 $ is spanned by the vector $v_{k+1}$. Let
$v^0(V^\bot_1) \in V^\bot_1 \cap V^\bot_2  $ and $v^0(V^\bot_2) \in V^\bot_1 \cap V^\bot_2 $ such that $\{ w_1+v_1, \dots, w_k+v_k, v^0(V^\bot_1)\}$
(resp.  $\{ w_1-v_1, \dots, w_k-v_k, v^0(V^\bot_2)\}$ ) is a  positively oriented basis of $V^\bot_1 $ (resp. $V^\bot_2$.) We have $ v^0(V^\bot_1) = (-1)^k v^0(V^\bot_2)$ since
$\pi_X(w_i+v_i)= v_i $ and $\pi_X(w_i - v_i)= - v_i $  for all  $i, i=1, \dots, k .$  In
particular,  $ v^0(V^\bot_1) = -v^0(V^\bot_2)$ when $k=1.$
The general case follows since any pair of maximal $B$--isotropic transversal subspaces of $\mathbb{R}^n$ is of the form $( g (V_1), g (V_2) )$ for some $g \in SO(B).$\\
\noindent\textbf{6.4 Margulis's sign.} Let us recall now the definition of the Margulis sign (or for short \textit{sign})  of an affine element. Let $g \in G_n$ be an $\mathbb{R}$-regular element with $l(g) \in SO(B)$ where $B $ is a non-degenerate form  on $\mathbb{R}^n$ of signature $(k+1,k).$ Obviously, the subspaces $A^+(g)$ and $A^-(g)$ are maximal $B$--isotropic subspaces, $D^+(g)= A^+(g)^\bot$ and
$D^-(g)= A^-(g)^\bot.$ Following the procedure above for the element $g$ we choose and fix  a vector $ v_+ $ with the following property $v_+=v^0(D^+(g)).$ 
 We did not  choose an orientation on the line  $A^0(g)$ at this point, but we fix an orientation on this line by the choice of the orientation on $A^+(g)$ and $D^+(g).$ Thus we fix an orientation on $A^0(g^{-1})$  by the choice of the orientations on $A^+(g^{-1})= A^-(g)$ and $D^+(g^{-1})= D^-(g).$  We will denote the corresponding vector by $v_-.$ Let $q \in \mathbb{R}^n.$  Set
$$\alpha(g) =B(gq - q, {v_+}) / B({v_+},{v_+})^{1/2}.$$ It is clear that $\alpha(g)$ does not depend on the point $q \in \mathbb{R}^n$ and  we have $\alpha(g) =\alpha (x^{-1} g  x)$
 for every $x \in G_n$ such that $l(x) \in SO(B).$  Consider now any $\mathbb{R}$-regular element $g$ and let us show that  $\alpha (g^{-1})= (-1)^k \alpha(g).$ Indeed by Example 6.3, $v_-=v^0(D^+(g^{-1})) =v^0(D^-(g))= (-1)^k v^0(D^+(g)) =(-1)^k v_+$. 
 We have $ \alpha(g^{-1})= B(g^{-1}q - q, {v_-}) / B({v_-},{v_-})^{1/2}= (-1)^k B(g^{-1}q - q, {v_+}) / B({v_+},{v_+})^{1/2} =
(-1)^{k+1}B(q - g^{-1}q, {v_+}) / B({v_+},{v_+})^{1/2}.$ Put $p =g^{-1}q.$  Hence $\alpha (g^{-1})= (-1)^{k+1}\alpha(g).$
\begin{figure}[!h]
    \label{pic-gamma}
 \centerline{ 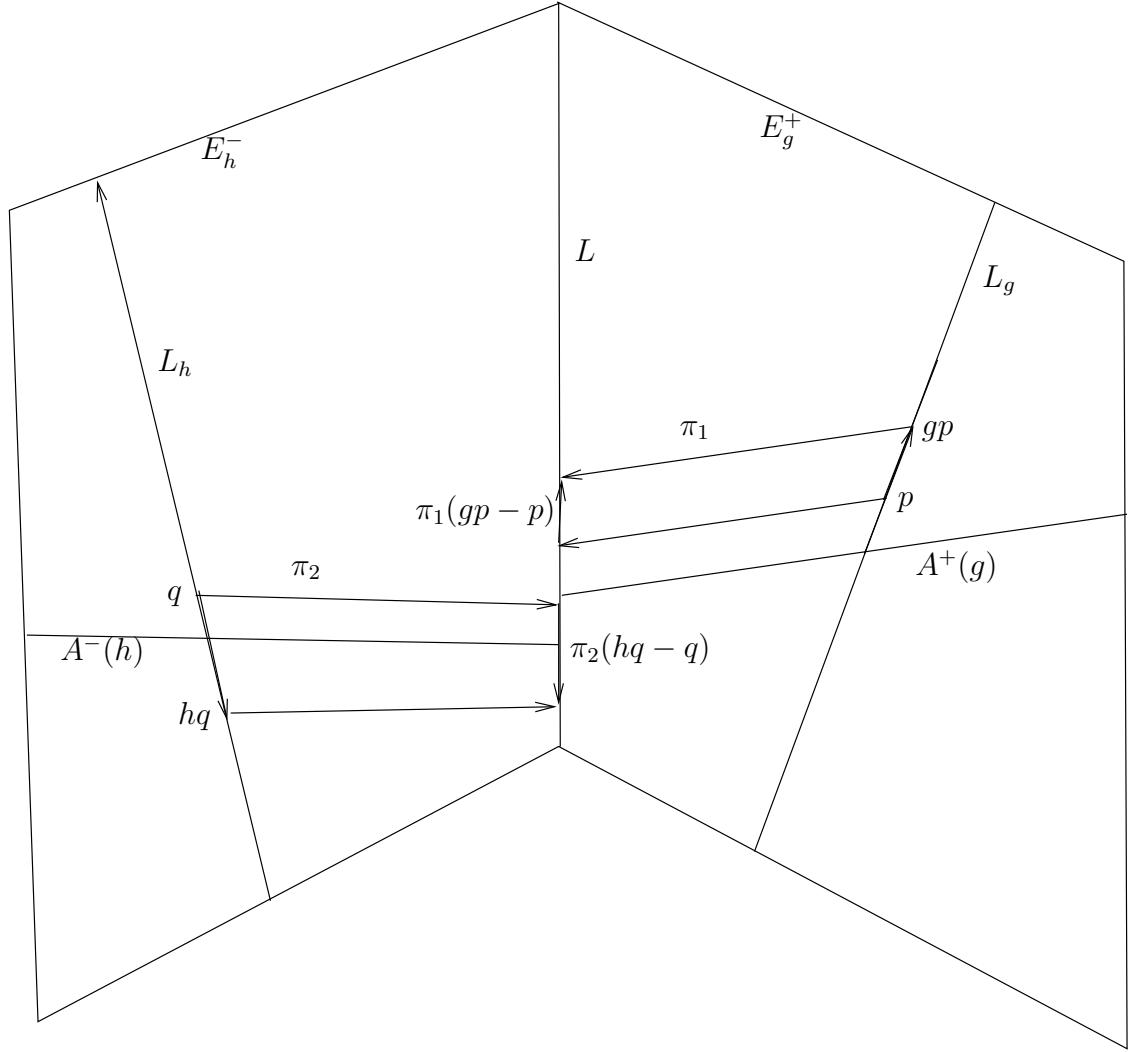} \vspace{15pt}
    \caption {Positive and negative parts, illustration 1}
   \end{figure} 
Note that  $\alpha(g) =\alpha(g^{-1}) $ if $k =1$.The number $\alpha(g)$ is called the sign of $g$ (see [AMS3]).
We call $\alpha(g)$ the sign of $g$, although $\alpha(g)$ is a non-zero real number, since we are only interested in its sign and not in its absolute value.
\par Using this approach we define now the sign of a regular element $g$ of the group $\Gamma$ for the case that the semisimple part of the Zariski closure $G$ of $\Gamma$ is $SO(2,1) \times SL_3(\mathbb{R}).$  Recall that $V =V_1 \oplus V_2$, $ S |_{ V_1} = SO(2,1)$ and $S |_{ V_2}= SL_3(\mathbb{R}).$ Hence we have two natural homomorphisms: $\theta_1 : G \rightarrow SO(2,1) \subseteq GL(V_1) $ and   $\theta_2 : G \rightarrow SL_3(\mathbb{R}) \subseteq GL(V_2).$ We will also assume that our standard inner product (see 2.4)\, is chosen so that the subspaces $V_1$ and $V_2$ are orthogonal.  
As the first step we have to choose the positive vector $v_g,  v_g \in A^0(g).$ Let $g \in S$ be a regular element. Let $\widehat{g}$ be the restriction $g| _{V_1}.$ Obviously $A^0(g) =A^0(\widehat{g}).$ We choose and fix a vector $v_+ =v^0(D^+(\widehat{g}))$ as we explained above. Set $v_g =v_+/B(v_+, v_+)^{1/2}.$ Let $g\in G$ be a regular element. 
There exists a unique $u \in G$ such that $l(h) =l(ugu^{-1}) \in S.$ Set $v_h =l(u) (v_g) .$ There is a simple geometrical explanation of this definition. Let $\pi : V \longrightarrow V_1$ be the natural projection onto $V_1$ along $V_2$. We have the corresponding homomorphism  $\widehat{\pi} : G \longrightarrow SO(2,1).$ It is easy to see that the restriction of $\pi$ to $A^0(g)$ gives an isomorphism onto $A^0(\widehat{\pi}(g))$ and $\pi (v_g) = v_{\widehat{\pi}(g)}.$
Let $\tau_g : V  \longrightarrow L_g  $ be the natural projection of the affine space $V$ onto
the line $L_g$ along the subspace $A^+(g) \oplus A^-(g)$, where $g $ is a regular element . There exists a unique $\alpha \in \mathbb{R}$ such that $\tau_g(p)-p = \alpha v_g.$ We set $\alpha (g) =\alpha.$  Clearly $\alpha (g) = B( \pi(\tau_g(p)-p) , \pi(v_g) )$ where $B$ is the form of signature $(2,1)$ on
$V_1$ fixed by $SO(2,1)$ since  $\pi(v_g) = v_{\widehat{\pi}(g)}.$  Obviously  $\alpha (g) $ does not depend on the chosen  point  therefore we have $\alpha (g^{-1})=\alpha (g)$  and $\alpha (g^n) =|n| \alpha (g) .$  For more details see [AMS4, p.5]. A regular element $g \in G$ is called hyperbolic, if $\theta_1(g)$ and $\theta_2(g)$ are hyperbolic.
\par Let us now explain the main application of these definitions. Let $g$ and $h$ be two hyperbolic  transversal elements.
Then $A^-(h) \oplus D^+(g)= V$ and $dim (D^-(h) \cap D^+(g)) =1.$ Let $A =D^-(h) \cap D^+(g)$. Let $L$ be the corresponding line $L = E^+_g \cap E^-_h.$ Let $\pi_1 :L_g \longrightarrow  L $ be the projection of $L_g$ onto $L$ along $A^+(g)$ and let    $\pi_2 : L_h \longrightarrow L. $ be the projection of $L_h$ onto $L$ along $A^-(h)$ (see Fig.2).   By the above arguments for  $p \in L_g, q \in L_h$ the vectors  $\pi_2(hq -q)$ and $\pi_1(gp -p)$ have opposite directions if $\alpha(g)\alpha(h) <0.$
Then as in the proof of Theorem A [AMS3], we conclude that there exist infinitely many positive numbers $n, m$ and two  balls $B(p,1)$ and $B(q,1)$ such that $h^mg^n B(p,1) \bigcap B(q,1)\neq \emptyset.$  Thus we conclude\\
\noindent \pro {\it Lemma  6.5}.  Let $\Gamma$ be an affine group and let $G$ be the Zariski closure of $\Gamma.$  Assume that the semisimple part $S$ of $G$ is as in the  Case 2 (3), and there are two hyperbolic transversal elements $g, h \in \Gamma$ such that $\alpha(g)\alpha(h) < 0 .$ Then the group $\Gamma$ does not act properly discontinuously.\endpro \\
\textbf {6.6}  To construct transversal hyperbolic elements of the group $\G$ with opposite sign is more difficult here than in the case when the semisimple part is $SO(2,1)$ (see [S1]). To make products transversal, one needs a quantitative version of hyperbolicity and transversality, see Lemma 2.7.  Thus we construct the appropriate set $M \subseteq \G$ of hyperbolic elements  to insure that a given hyperbolic element $\g \in \G$ will be at least $\varepsilon =\varepsilon(M)$--transversal 
to some element of $M.$ 
 Moreover, the set $M$ will be the union of two sets $M_1$ and $M_2.$ 
If  the number of eigenvalues of $\g$ greater  than $1$ is two (resp. three) then  $\g \in \G$ will be $\varepsilon$--transversal  to at least one  element of $M_1$ (resp. $M_2$). This is close to the strategy we used in [AMS1].
\par Recall that $ v_1, v_2, w_1$ is a basis of $V_1$ such that for any vector $v \in V_1, v = x_1v_1 +x_2v_2 +y_1w_1$ we have $B(v,v) = x_1^2 +x_2^2 -y_1^2.$ 
We will use the notations and definitions from 6.1. Let $\partial {Con_B}$ be the  boundary of  $Con_B.$ 
\begin{figure}[!h] 
    \label{pic-gamma}
  \centerline{ 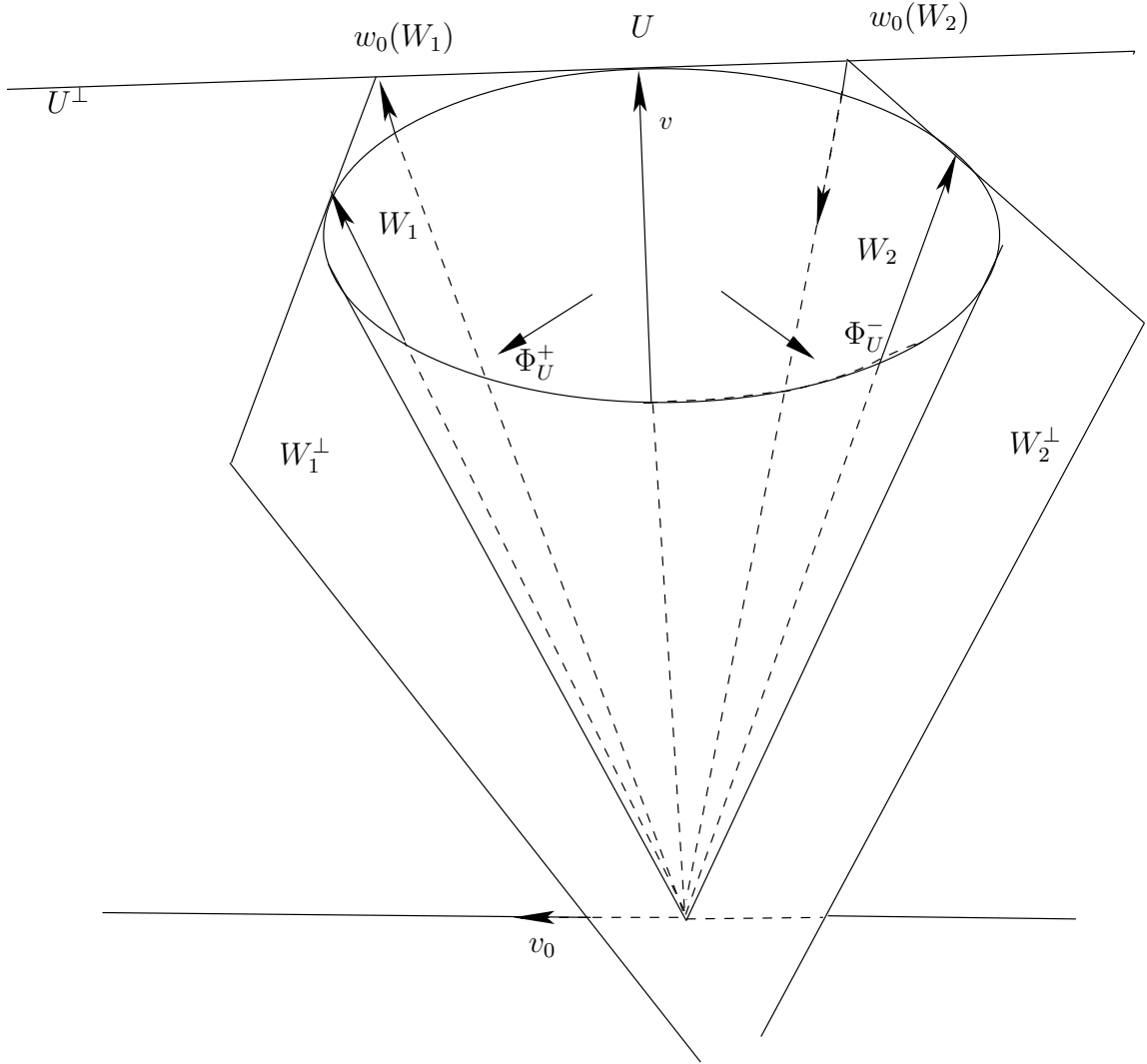} \vspace{15pt}
   \caption {Positive and negative parts, illustration 2}
  \end{figure}
Let $U$ (see Fig.3)  
  be a maximal $B$-isotropic subspace of $V_1$ and let $v$ be the vector of $U$ such that  $\pi_Y(v)= w_1 .$ Clearly  $U$ is spanned by $v.$   Let $v_0$ be the vector in $U^\bot \cap X$ such that $B(v_0,v_0)=1$ and the basis $\pi_X (v), v_0$ has the same orientation as $v_1,v_2$.  
Let $W$ be a maximal $B$-isotropic subspace of $V_1$ and suppose $W \neq U$ . Then  $ \dim (U^ \bot \cap W^ \bot)=1$. There exists a unique vector $w_0(W)$ in $ U^\bot \cap W^ \bot$ and $\widehat{v} \in U$ such that $w_0(W) = v_0 + \widehat{v}$. Obviously there exists a unique number $\alpha(W)$ such that $\widehat{v} =\alpha(W) v.$  Set  $\Phi^+_U = \{W \in \Phi |\, \alpha(W) > 0   \}$ and  $\Phi^-_U = \{W \in \Phi |\, \alpha(W) < 0   \}.$ We have $B(v_0,w_1)=0$ since $v_0 \in X$. Therefore $B(w_0(W)),w_1) =\alpha(W) B(v,w_1) = -\alpha(W)$. Let $\widehat{U}$  be the sum of the two subspaces $U$ and $< w_1>$. Then $\Phi^+_U $ and $\Phi^-_U$ are two different connected components of the set
$\partial {Con}_B\backslash \widehat{U}.$ Obviously $\partial {Con}_B\backslash\widehat{U} = \Phi^+_U \cup \Phi^-_U. $
We conclude :  \\
$(1)$ For every $W \in \Phi^+_{U} $ (resp. $W \in \Phi^-_{U} $) we have $B(w_0(W), w_1) <  0$  (resp. $B(w_0(W) ,w_1) >0$).\\
$(2)$ Let $W_1, W_2, W_3, W_4$ be maximal $B$--isotropic subspaces  of $V_1$  such that  $ w_1 \in (W_1 +W_2) \bigcap (W_3+W_4).$ It is easy to see that if $U \notin \{ W_1, W_2\}$ then $W_1$ and $W_2$ belong to different connected components
of the set  $\partial {Con}_B \setminus \widehat{U}$. Indeed, since $w_1 \in W_1+W_2$ we have $\alpha(W_1) = -\alpha (W_2).$ The same is true in the case $U \notin \{ W_3, W_4\}.$\\
$(2_a)$ It follows from (2) that  for every maximal $B$--isotropic subspace $U$ of $V_1$ if $W_{i} \in \Phi^{\pm}_ {U}$ then  $W_{i+1} \in \Phi^{\mp}_ {U}$ where $i=1,3.$\\
\noindent $(2_b)$  Let $d = \min _{1 \leq i \neq j \leq 4}\{d(W_i, W_j)\}.$ Let  $U$ be a maximal  $B$--isotropic subspace of $ V.$ It follows from $(2_a)$ that there exists $\delta =\delta(d)$ such that for every four maximal  $B$--isotropic subspaces $\widehat{W}_i, i=1,2,3,4$ of $V$ with $d(\widehat {W}_i, W_i) \leq \delta$ for  $ 1\leq i\leq 4$  there exists an $i_0 \in\{1,3\}$ such that $\widehat{W}_{i_0} \in \Phi^-_ {U} $ and \\$\widehat{W}_{i_0+1} \in \Phi^+_ {U} .$ \\
$(3)$ Assume first that $W_1 \in \Phi^+_{U} $ and    $W_2 \in \Phi^-_{U}.$ Let $\varepsilon = \min \{d( W_1, U), d (W_2, U)\}. $ Then there exists a $\delta =\delta(\varepsilon)$ such that if $\widehat{U}$ is a maximal $B$--isotropic subspace with $d(\widehat{U}, U) < \delta$ we have
$ \widehat{W}_1 \in \Phi^+_ {\widehat{U}}$ and $ \widehat{W}_2 \in \Phi^-_ {\widehat{U}}$ for any maximal $B$--isotropic subspaces  $\widehat{W}_1,  \widehat{W}_2$  with $d(\widehat{W}_1, W_1) < \delta$ and $d(\widehat{W}_2, W_2) < \delta.$\\ 
\pro{\it Remark 3.}\endpro Let us explain the  motivation of  the above construction. Assume that we have two hyperbolic transversal elements $\gamma_1$ and $\gamma_2$ of $\Gamma$ such that 
\begin{enumerate}
\item[(1)] $A^+(\theta_1(\gamma_1))= A^+(\theta_1(\gamma_2)) = U$
\item[(2)] $A^-(\theta_1(\gamma_1)) \in \Phi^+_U$ and  $A^-(\theta_1(\gamma_2)) \in \Phi^-_U$ 
\item[(3)] There exists a point $x_0 \in \mathbb{R}^6$ such that the vectors  
$$v_i =\frac{\gamma_i x_0 - x_0}{d(\gamma_i x_0 , x_0 )} $$ 
\end{enumerate}
are positive multiples of $w_1, i =1,2.$\\ Remember that $w_1$ is the third basis vector of our standard basis of $V_1.$
We conclude from (1), (2) and (3) that $$\alpha(\gamma_1) \alpha (\gamma_2) = B(v^0_{ \theta_1(\gamma_1)}, v_1)B(v^0_{ \theta_1(\gamma_2)}, v_2) < 0.$$
Thus the two hyperbolic elements $\gamma_1$ and $ \g_2 $ have opposite sign.
\par Note that the conclusion about the sign is an open condition. Thus if we vary  $\gamma_1$ and $\gamma_2$ slightly then this conclusion remains true. In fact, there exists a positive number $\delta$ such that if  $\gamma_1$ and $\gamma_2$  are two hyperbolic transversal elements of $\Gamma$ such that
\begin{enumerate}
\item[(1)] $d(A^+(\theta_1(\gamma_1)),U) < \delta$ and $d(A^+(\theta_1(\gamma_2)), U) < \delta$ for $U$ as above
\item[(2)] $A^-(\theta_1(\gamma_1)) \in \Phi^+_U$ and  $A^-(\theta_1(\gamma_2)) \in \Phi^-_U$  and
\item[(3)] there exists a point $x_0 \in \mathbb{R}^6$ such that for $i=1,2$ the angle $\angle (v_i,w_1) < \delta,$ where $$v_i =\frac{\gamma_i x_0 - x_0}{d(\gamma_i x_0 , x_0 )}$$ 
then $\alpha(\gamma_1) \alpha(\g_2) <0.$
\end{enumerate}
\noindent \pro{\it Lemma 6.7} Let $\widehat{\Gamma} \subset GL(V_1)$ be a Zariski dense subgroup of $SO(2,1).$ Then there exist four transversal hyperbolic elements $\g_1, \g_2,\g_3, \g_4$ in $\widehat{\Gamma}$ such that 
we have $B(v,v) <0$  for every non- zero vector 
$$v \in (A^+(\g_1) +A^+(\g_2))\cap (A^+(\g_3) +A^+(\g_4))$$

\noindent \pro{\it Proof} \endpro Since $\widehat{\Gamma} $ is Zariski dense in $SO(2,1)$ there are four transversal hyperbolic elements [AMS1]. It is enough now to give an order of these four elements such that any non-zero vector  $v \in (A^+(\g_1) +A^+(\g_2))\bigcap (A^+(\g_3) +A^+(\g_4))$ will be inside the cone $Con_{B}.$ 
Thus $B(v,v) < 0$ for any non-zero  $v \in (A^+(\g_1) +A^+(\g_2))\bigcap (A^+(\g_3) +A^+(\g_4))$ which proves the lemma.
\par Let us show that there exist four transversal hyperbolic elements $\g_1, \g_2, \g_3$,\\$ \g_4$ of $\G$ such that $ w_1 \in (A^+(\g_1) +A^+(\g_2))\bigcap (A^+(\g_3) +A^+(\g_4)).$  
\endpro
\begin{figure}[!h] 
    \label{pic-gamma}
\centerline{ \input{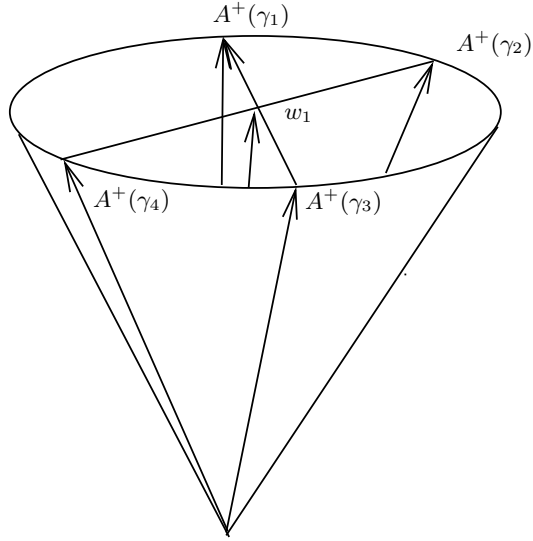}} \vspace{15pt}
    \caption {Configuration }
   \end{figure}
Indeed, let  $\widehat{\G} =\theta_1(\G).$ Recall that we are dealing with case 2 (3). It follows from Lemma 6.7 that there exist four transversal hyperbolic elements $\g_1, \g_2,\g_3, \g_4$ of $\G$ such that $v \in (A^+(\theta_1(\g_1)) +A^+(\theta_1(\g_2)))\bigcap (A^+(\theta_1(\g_3)) +A^+\theta_1((\g_4)))$, $B(v,v)=-1.$ Obviously, there exists an element $g \in O(2,1)$ such that $g v =w_1.$ Let  $\overline{\G} =g \G g^{-1}$ and let $\overline{\g}_i =g \g_i g^{-1}, i=1,2,3,4.$ We conclude that  $ w_1 \in (A^+(\overline{\g}_1) +A^+(\overline{\g}_2))\bigcap (A^+(\overline{\g}_3) +A^+(\overline{\g}_4)).$ Using the same notation we can and will assume that \\
  $w_1 \in (A^+(\theta_1(\g_1)) +A^+(\theta_1(\g_2)))\bigcap (A^+(\theta_1(\g_3)) +A^+(\theta_1(\g_4)))$ (see Fig.4 above).  Set $A_i =A^+(\theta_1(\g_i)),$ for $ i=1,2,3,4.$
\par Let us explain in more detail how we prove that in Case 2 (3) $\Gamma$ is not crystallographic. Assume that $\G$ is crystallographic and let $K$ be a compact subset of $\mathbb{R}^6 $ such that $\G K = \mathbb{R}^6 .$ \\
 Let $S \subseteq$ \text{Aff}$(\mathbb{R}^n)$ be an infinite subset. We will say that $S$ is unbounded if for every compact $K, K \subset \mathbb{R}^n$  there exists an infinite subset $\tilde{S}$ of $S$ such that $s K \cap K = \emptyset$ for all $s \in \tilde{S}.$ Every infinite subset  $S$ of $\G$ is unbounded since $\G$ acts properly discontinuously.
\par  We will say that a non-zero vector $v \in \mathbb{R}^6 $ is the direction of an unbounded subset  $S$ on a compact $K_0$ if there exist a positive real number $t$, an unbounded sequence $\{\gamma_n, n\in \mathbb{N}\}$ of $S $, and a sequence $\{ x_n , n \in \mathbb{N} \} $ in $K_0$ such that  $$ \frac{\gamma_n x_n -x_n}{{d(\gamma_{n}x_n , x_n )} } \rightarrow t v $$ for $n \rightarrow \infty$. Denote by $V(S, K_0 )$ the set of all directions of $S$ on $K_0$. 
\par Every non zero vector in $\mathbb{R}^6$ is contained in $V(\G,K)$ since $\G$ is crystallographic.\\ 
Let $\g \in \G$ and $\tilde{S}=S \g.$ Set $\tilde{K}= \g^{-1}K.$ Then 
$V(S, K )= V(\tilde{S}, \tilde{K} ).$ Let $M= \{ \g^{(1)}, \dots \g^{(m)}\}$ be a finite subset of $\G.$ Set $\widehat{S}= \cup ^{m}_1 S\g^{(i)} $ for $i=1, \dots m.$ 
Let $\widehat{K} =\cup ^{m}_1 (\g^{(i)})^{-1} K.$ Then $V(S, K ) \subseteq V( \widehat{S}, \widehat{K} ).$ \\
We will show (see the proof of the main lemma) that the "finiteness theorem" [AMS1] and  the  arguments above imply that a given unbounded set $S \subset \G$  can be transformed into an unbounded  set of $\varepsilon$-hyperbolic elements of $\G$, where $ \varepsilon= \varepsilon(\G),$  such that for every positive $\delta$ the set $\{\g \in S, s({\gamma})> \delta\} $ is finite. At the same time we can transform the compact subset $K$ into a compact subset  $\widehat{K}$ such that $V (S, K) \subseteq V(\widehat{S}, \widehat{K}).$ \par Let   $S$  be a subset
 such that $w_1 \in  V(S,K) .$  We thus  can and will assume that $S$  consists of $\varepsilon$-hyperbolic elements such that for every infinite subsequence $\{\gamma_n\}$ of $S$ we have
 $s({\gamma}_n)\rightarrow 0$ for $n \rightarrow \infty.$\\
   The sequence of isotropic one dimensional spaces $A^+(\theta_1(\gamma_n))$ converges to an isotropic subspace $W$ since the projective space is compact, Then by Lemma 6.8 below we will conclude that there exist at least two elements $\gamma_{+}$,
 $\gamma_{-}$  and a positive number $M \in \mathbb{Z}$ such that for $m > M$ and all $n$ we have
\begin{enumerate}
\item[(1)] $A^-(\theta_1(\gamma_n \gamma^m_{+})) \in \Phi^+(W)$ and $A^-(\theta_1(\gamma_n \gamma^m_{-})) \in \Phi^-(W)$ 
\item[(2)]  $\lim_{n \rightarrow \infty} A^+(\theta_1(\gamma_n \gamma^m_{\pm}) )= W$  for every $m > M$.
\end{enumerate} 
Fix $m_0 >M$.  The arguments above imply that $ V(\{\gamma_n \g^{m_0}_{\pm}\}, \g^{-m_0}_{\pm} K )$ contains $w_1.$  
 It follows  from Remark 3 and  (1) and (2) above that for sufficiently big $n$ we have $\alpha(\gamma_n \gamma^{m_0}_{+})  \alpha(\gamma_n \gamma^{m_0}_{-}) < 0. $ Hence there are two elements in $\G$ with opposite sign. We thus can conclude that $\Gamma$ does not act properly discontinuously. \\
 Let us end with a small remark. We explained above how to change a sign and how to obtain two elements with different signs. On the other hand, to be able to use  Lemma 6.5 we need two hyperbolic elements. The difficulty which comes up here is the following. To conclude that the product $\g_n \g^m_{\pm}$ is hyperbolic we need transversal elements  $\g_n$ and $\g^m_{\pm}.$ For $\g \in \G$ we have $\theta_2(\g) \in SL_3(\mathbb{R}).$ Therefore  $\dim A^-(\theta_2(\g_n)) =2 $ or  $\dim A^-(\theta_2(\g_n)) =1.$ Since we do not know a priory the dimension of  $\dim A^-(\theta_2(\g_n)) $ we "prepare" two sets $S_i$ and $T_i$ which fulfil (4) of Lemma 6.8  below.\\
\noindent \pro{\it Lemma 6.8} For every subspace $A_i , 1 \leq  i \leq 4,$ and positive $\delta$ there exist sets $S_i =\{g_{i1}, g_{i2}, g_{i3} \} \subset \Gamma$ and $T_i =\{h_{i1}, h_{i2}, h_{i3}\} \subset \Gamma$  and positive real numbers  $\varepsilon, q <1$, such that
\begin{enumerate}
\item[(1)] $\widehat{d}(A^-(\theta_1(g_{ik})), A_i) < \delta $ and  $\widehat{d}(A^-(\theta_1(h_{ik})), A_i) < \delta; $
\item[(2)] $g_{ik}$ and $h_{ik}$ are $\varepsilon$-hyperbolic for $k =1,2,3;$
\item[(3)] $\max_{1\leq i\leq 4, 1 \leq k \leq 3} \{s(g_{ik}), s(h_{ik})\} <q;$
\item[(4)] let $i$ be an index with $1\leq i \leq 4.$ Then for every $k=1,2,3$ we have\\ $\dim A^-(\theta_2(g_{ik})) =2$ and $\dim A^-(\theta_2(h_{ik})) =1;$
\item[(5)] for every index $i$ with $1\leq i \leq 4$ we have $\bigcap_{1 \leq k \leq 3} A^-(\theta_2(g_{ik})) = \{0\};$
\item[(6)] for every index $i$ with $1\leq i \leq 4$ we have $\dim (A^-(\theta_2(h_{i1})) +A^-(\theta_2(h_{i2}))+ A^-(\theta_2(h_{i3})) ) =3. $
\end{enumerate}
\endpro
\noindent \pro{\it Proof} \endpro Obviously  it is enough to prove the statement for one subspace. Let us do it for $A_1.$
It is easy to show that there exists a hyperbolic element  $\gamma$ of $\Gamma$ such that 
\begin{enumerate}
\item[i)] $\theta_1(\gamma)$ and $\theta_1(\gamma^{-1})$ are transversal to $\theta_1(\gamma_1);$
\item[ii)] any  proper $\theta_2(\gamma)$--invariant subspace  does not contain  a proper $\theta_2(\gamma_1)$--invariant subspace;
\item[iii)] any proper $\theta_2(\gamma_1)$-- invariant subspace does not contain  a proper $\theta_2(\gamma)$--invariant subspace.
\end{enumerate}
We will also assume that   $\theta_2(\gamma)$ has three eigenvalues of different norms  [AMS1]. In that case all of them are real numbers.
Put $\gamma_n = \gamma^n_1 \gamma \gamma_1^{-n}.$ We can assume that $\dim A^-(\gamma) = 2$ otherwise we can take
$\gamma^{-1}$ instead of $\gamma.$ Let us first show that for some positive numbers $n_1, n_2, n_3$ we have
$ \cap_{1 \leq i \leq 3}  A^-(\theta_2(\gamma_{n_i})) =\{0\}.$ Since  for $n \neq m$ we have $ A^-(\theta_2(\gamma_{n})) \neq  A^-(\theta_2(\gamma_{m}))$ there
are positive numbers $n_1$ and $n_2$ such that $\dim A^-(\theta_2(\gamma_{n_1})) \cap A^-(\theta_2(\gamma_{n_2}) =1 .$ Let $v$ be a non -zero vector of this intersection. If
$\theta_2(\gamma_1)^{-n} v \in A^-(\theta_2(\gamma))$ for infinitely many positive $n $ then the proper $\theta_2(\gamma)$--invariant subspace $ A^-(\theta_2(\gamma))$ contains a  $\theta_2(\gamma_1)$-- invariant subspace. This contradicts our assumptions. Thus, by the choice of $\gamma$ and $\gamma_1$ there 
exists an $n_3$ such that  $\theta_2(\gamma_1)^{-n_3} v \notin A^-(\theta_2(\gamma))$
 Therefore $v \notin\theta_2(\gamma_1)^{n_3} A^-(\theta_2(\gamma)) =  A^-(\theta_2(\gamma_{n_3})) .$ Clearly\\
 $ A^-(\theta_2(\gamma_{n_1+m})) \cap $ $A^-(\theta_2(\gamma_{n_2+m})) \cap$ $A^-(\theta_2(\gamma_{n_3+m})) =\{0\}  $ for all positive numbers $m.$ Since the projective space $PV$ is compact we can and will assume that $A^+(\gamma_{n_i+m}) \longrightarrow X_i^+$,  $A^-(\gamma_{n_i+m})\longrightarrow X^-_i$ for $m \longrightarrow \infty$ and $i=1,2,3.$ 
  By standard arguments [MS], [AMS1], we conclude that there exists a hyperbolic element $\gamma_0$ such that
 $\widehat{d} (A^+(\gamma_0), X^-_i) > 0$ and $\widehat{d} (A^-(\gamma_0), X^+_i) > 0$ for all $i =1,2,3.$ \\Let
 $\widehat{\varepsilon} = \min _{1 \leq i \leq 3}\{ \widehat{d} (A^+(\gamma_0), X^-_i) , \widehat{d} (A^-(\gamma_0), X^+_i) \}.$
 Thus there exists an $M \in \mathbb{N}$ such that for $m \geq M$ the elements $\gamma_0$ and $\gamma_{n_i+m} $
 are  $\widehat{\varepsilon}/2 $--transversal.
 Let $q_1 =\max \{s(\gamma_0),$  $ s(\gamma)\} .$ 
Since $ \g_m =\g^m_1\g \g_1^{-m}$ we have $ A^+(\theta_1( \g_m)) = \theta_1(\g_1)^m A^+(\theta_1(\g)).$ Therefore 
the sequence $\{ A^+(\theta_1(\g_m) \}_{m \in \mathbb{N}}$ converges to $A^+(\theta(\g_1)) =A_1$ when $m \rightarrow \infty$
(see for example [AMS1]). Hence for every positive $\delta$ there exists an $M$ such that  for $m > M$ we have $\widehat{\rho}(A_1, A^-(\theta_1(\gamma_m))) \leq \delta/4.$ 
It follows from [MS] and [AMS1] that there exists a positive number $N$ such that for $n >N$ we have $$\widehat{\rho}(A^+(\gamma^n_0\gamma^n_i, A^+(\gamma_0)) \leq q_1^{n/2}$$ $$\widehat{\rho}(A^-(\gamma^n_0\gamma^n_i), A^-(\gamma_i)) \leq q_1^{n/2}$$ and
$s(\gamma^n_0\gamma^n_i) \leq q_1^{n/2}$ for $i =1,2,3.$ Therefore there exist positive numbers $N_1$ and $\varepsilon$ such that for all $i =1,2,3$ the element $\gamma^n_0\gamma^n_i$  is $\varepsilon/2$-hyperbolic
and $ \widehat{\rho}( A_1, A^+(\theta_1(\gamma^n_0\gamma^n_i)) < \delta/2$ for $n > N_1.$ Since $ A^-(\gamma^n_0\gamma^n_i) \longrightarrow A^-(\gamma_i)$ for $n \longrightarrow \infty$, there exists an $N \in \mathbb{N}$ such that  $\dim A^-(\theta_2(\gamma^n_0\gamma^n_i))  =$ $\dim  A^-(\theta_2(\gamma^n_i) =2$ for $ i=1,2,3$ and  for $n >N_2$  we have $\cap_{1 \leq i \leq 3} A^-(\theta_2(\gamma^n_0\gamma^n_i)) = \{0\}.$  Take $n > \max\{N_1, N_2\}$ and set $g_{1k} = \gamma^n_0\gamma^n_k , k=1,2,3$. Using the same arguments one can show that there is a set  $T_1 =\{h_{11}, h_{12}, h_{13}\} $ with properties 1-4, 6. This proves Lemma 6.8.
\par For chosen sets $S_i =\{g_{i1}, g_{i2}, g_{i3} \} \subset \Gamma$ and $T_i =\{h_{i1}, h_{i2}, h_{i3}\} \subset \Gamma , i=1,2,3,4$ we will define the following constants. For any one dimensional subspace $U$ of $V_1$, we have
$\sum_{1 \leq i \leq 4, 1 \leq k \leq 3}\widehat{d}(U, A^- (\theta_1(g_{ik}))) > 0. $ Since the projective space is compact we have $$\inf_{U \subset V_1} \sum_{1 \leq i \leq 4, 1 \leq k \leq 3}\widehat{d}(U, A^-(\theta_1(g_{ik})))  >0.$$
Set $$d_1^{(S)} =\inf_{U \subset V_1} \sum_{1 \leq i \leq 4, 1 \leq k \leq 3}\widehat{d}(U, (A^-(\theta_1(g_{ik}))) /100.$$ By the same arguments there exists a positive constant $d_1^{(T)}$, such that $$d_1^{(T)} = \inf_{U \subset V_1} \sum_{1 \leq i \leq 4, 1 \leq k \leq 3}\widehat{d}(U, A^-\theta_1((h_{ik}))) /100.$$
\par Let $U$ be a one dimensional subspace of $V_2$. It follows from Lemma 6.8, (5)  that $\sum_{1 \leq i \leq 4, 1 \leq k \leq 3}\widehat{d}(U, A^-(\theta_2(g_{ik}))) > 0. $ Thus 
$$\inf _{U \subset V_2 }\sum_{1 \leq i \leq 4, 1 \leq k \leq 3}\widehat{d}(U,A^-( \theta_2(g_{ik}))) > 0.$$  Set $$d_2^{(S)} = \inf _{U \subset V_2 }\sum_{1 \leq i \leq 4, 1 \leq k \leq 3}\widehat{d}(U, A^-(\theta_2(g_{ik}))) /100.$$ Let $U$ be a two-dimensional subspace of $V_2$. It follows from  Lemma 6.8, (6)  that  $\sum_{1 \leq i \leq 4, 1 \leq k \leq 3}\widehat{d}(U, A^-(\theta_2(h_{ik}))) > 0. $ Now by the same arguments as above there exists a positive
$$d_2^{(T)} = \inf _{U \subset V_2} \sum_{1 \leq i \leq 4, 1 \leq k \leq 3}\widehat{d}(U, A^-(\theta_2(h_{ik})))/100,$$
\noindent\pro {\it Main Lemma 6.9 } There are two hyperbolic elements of the group $\Gamma$ such that  $\alpha (g) \alpha (h) < 0. $ \endpro\\
\pro {\it Proof.} \endpro We have to prove that there are two elements with opposite sign. Recall that we are dealing with case 2 (3). Since we can and will assume that there exists a hyperbolic element of positive sign, we will prove that there exists an element with negative sign.
\par Let $S_i =\{g_{i1}, g_{i2}, g_{i3} \} \subset \Gamma,$ $T_i =\{h_{i1}, h_{i2}, h_{i3}\} \subset \Gamma$ and  positive real numbers $\varepsilon, q <1 $ be as in Lemma 6.8. Let $d = \min _{1 \leq i \neq j \leq 4}\{\widehat{d}(A_i, A_j)\}$ and let $\delta (d)$ be as in $(2_b).$   Assume that we choose  a positive $\delta$ in  Lemma 6.8, (1)  such that $\delta \leq \delta(d)/4$. Set $$\varepsilon_1 = \min \{ d^{(S)}_1, d_1^{(T)}, d_2^{(S)}, d_2^{(T)} \}.$$
\par Let $K$ be a compact subset of $V$ such that $\Gamma K = V.$ Denote by $L$ the ray $L= \{t w_1, t \in \mathbb{R}, t >0 \}.$ We may assume that $K \bigcap L \neq \emptyset.$ Then there exists a sequence  $\{\g_n\}_{n \in \mathbb{N}}$ of elements of $\G$ and a sequence of points $p_n \in L$ such that
\begin{enumerate} 
\item[(1)] $\g^{-1}_np_n \in K$
 \item[(2)] $d(p_n,\g^{-1}_np_n ) \rightarrow \infty$ when $n \rightarrow \infty.$
 \end{enumerate}
Set $k_n = \g^{-1}_np_n  \in K.$ It is easy to see that for $n \rightarrow \infty$ we have
\begin{enumerate} \item[(3)]
$\g_nk_n -k_n / d(\g_nk_n, k_n) \rightarrow w_1 .$
\end{enumerate}
 By [AMS1]  there exist an $\varepsilon _2=\varepsilon (\G )$ and a finite subset $S(\Gamma)=\{g_1, \dots, g_m\} $ of $\G$ such that  for every $\g \in \G$  there exists $g_i, i =i( \g), 1 \leq i \leq m$ and  $M =M(\varepsilon_2)$ such that the element $\g g^m_i$ is $\varepsilon_2$--hyperbolic and  $s(\g g^m_i) < s(g_i)^{m/2}$  for $m > N.$  We can choose an infinite subsequence $\g_{n_k}$ such that the element $g_{i} \in  S(\Gamma)$ is the same for all $\g_{n_ k}.$ Assume that this is $g_1.$
Put $r_m  =g_1^{-t}k_m $. This is a point of the compact set  $K_1 =g_1^{-t} K.$ Then for a fixed $t$ we have
 \begin{enumerate}
  \item[(4)] $g_1^{-t}\g^{-1}_np_n \in K_1$
 \item[(5)] $\g_n g^t_1 r_n -r_n / d(\g_n g^t_1 r_n , r_n) \rightarrow w_1$ for  $n \rightarrow \infty.$
 \item[(6)] $\lim_ {n \rightarrow \infty} s(\g_n g_1^t ) =0$
 \end{enumerate}
 Thus we assume that there exists a sequence  $\{\g_n\}_{n \in \mathbb{N}}$  of $\varepsilon_2$--hyperbolic elements of $\G,$ a compact set $K$ with  $K \cap L \neq \emptyset$ and a sequence of points $k_n \in K$ with properties (1),(2),(3) and (6).
The projective space $PV$ is compact. Thus we can and will assume that the sequences  $\{A^+(\g_n )\}_{n \in \mathbb{N}}$ and  $\{A^-(\g_n )\}_{n \in \mathbb{N}}$ converge. Let
$A^+(\g_n ) \longrightarrow A^+$ when $n \longrightarrow \infty$ and $A^-(\g_n ) \longrightarrow A^-$ when $n \longrightarrow \infty.$ Set $A_{\theta_1}^+ = V_1 \cap A^+, A_{\theta_1}^- = V_1 \cap A^-.$
\par There are two cases. 
\begin{enumerate}
\item[ (i)] The set $\{ n \in \mathbb{N} : \dim A^-(\theta_2(\g_n  ) =2 \}$ is infinite;
\item[(ii)] The set  $\{ n \in \mathbb{N}: \dim A^-(\theta_2(\g_n  ) =1 \} $ is infinite;
\end{enumerate}
 In case (I) we will consider the sets $S_i, i =1,2,3,4$. In case (II) we will consider the sets  $T_i, i =1,2,3,4$  and will use the following procedure. 
 \par Assume that  for infinitely many $n \in \mathbb{N}$ we have $\dim A^-(\theta_2(\g_n  ) )=2.$ It follows from [AMS1] that there exists a hyperbolic element $\gamma_0$ such that $\gamma_0$ and $g_{ik}$ are transversal for all $g_{ik } \in S _i, i =1,2,3,4,$\, $A^+(\gamma_0) \cap A^- =\{0\}$ and  $A^-(\gamma_0) \cap A^+ =\{0\}.$
Thus there exists $\varepsilon_3$ such that for every $n$ and every $g_{ik } \in S _i, i =1,2,3,4$ the two pairs of  hyperbolic elements $\{\g_n, \gamma_0\}$ and  $\{g_{ik }, \gamma_0\}$ are $\varepsilon_3$-transversal.  It follows from Lemma 2.7  that there exists a positive number $M \in \mathbb{N}$ such that for all $n \in \mathbb{N}$ the element $\g_n  \gamma^m_0$ is $\varepsilon_3/4$-hyperbolic and $\widehat{\rho}(A^-(\g_n  \gamma^m_0), A^-(\gamma_0)) \leq q^m_2 < \varepsilon/8$  and $s(\g_n  \gamma^m_0) \leq s(\gamma_0)^{m/2}$ for $m \geq M.$ Thus $\widehat{d}(A^-(\g_n \gamma^m_0)),A^+(g_{ik}) ) > \varepsilon_3/2 $ for all $g_{ik} \in S_i, i=1,2,3,4,$ $n \in \mathbb{N}$, $m \geq M.$ There exists an $M_1$ such that for $m \geq M_1$ we have
$ s(\gamma_0)^{m/2} \max \{\varepsilon_1,\varepsilon_3 \} \leq  \min  \{\varepsilon_1,\varepsilon_3 \} /8.$
 Fix $m > M_1,  m \in \mathbb{N}$ and set $\widehat{\gamma}_n =\g_n  \gamma^m_0.$
Obviously $\min_{n \in \mathbb{N}, 1 \leq i \leq 4, 1\leq k \leq 3} \widehat{d}(A^+(\theta_1( \widehat{\gamma}_n )), A^-(\theta_1 (g_{ik})) > \varepsilon_1$ and  $\min_{n \in \mathbb{N}, 1 \leq i \leq 4, 1\leq k \leq 3}
\widehat{d}(A^+(\theta_2( \widehat{\gamma}_n )), A^-(\theta_2(g_{ik})) > \varepsilon_1.$ It follows from 6.6 (2) that there exists an index $i_0, 1 \leq i_0 \leq 4,$ such that $ A^-(\theta_1 (g_{i_0 k})  \in \Phi^-_ {A^+(\theta_1( \widehat{\gamma}_n ))}.$ Without loss of generality, we will assume that $i_0 =1.$ Clearly $\min_{n \in \mathbb{N}, 1\leq k \leq 3} \widehat{d}(A^+(\theta_1( \widehat{\gamma}_n )), A^-(\theta_1 (g_{1k})) > \varepsilon_1/10.$ Then for some $k$ we have
$\widehat{d}(A^+(\theta_1( \widehat{\gamma}_n )), A^-(\theta_1 (g_{1k})) > \varepsilon_1/30.$ Assume this holds for $k=1.$  Thus $\widehat{d}(A^+(\theta_1( \widehat{\gamma}_n )), A^-(\theta_1 (g_{11})) > \varepsilon_1/30.$ On the other hand we know that\\
 $\widehat{d}(A^-(\widehat{\g}_n), A^+(g_{11}) ) > \varepsilon_3/2.$ It follows from (6), Lemma  2.7 (1)  and [MS], [AMS1]  that for any $m >0$ we have $$ A^+(\theta_1 (\widehat{\g}_ng_{11}^m)) \rightarrow A^+_{\theta_1} $$ when $n \rightarrow \infty.$
Therefore it follows from 6.6 (3) that there exists a positive number $N_0$ such that for $n > N_0$ and $\overline{\g}_n =\widehat{\g}_ng^{2N_0} _{11} $, $\varepsilon=
\min \{\varepsilon_1, \varepsilon_2\}/10$ we have
\begin{enumerate}
\item[(7)] $\overline{\g}_n $ are $\varepsilon$-hyperbolic elements ;
\item[(8)] $A^-(\theta_1(\overline{\g}_n)) \in \Phi _{A^+(\theta_1(\overline{\g}_n))}$;
\item[(9)] There exist a compact set $K_0$ and a sequence of points $\{k_n\} _{n \in \mathbb{N}} \subset K_0$
such that $\overline{\g}_n (k_n) \in L$ and $d ( \overline{\g}_n (k_n) , k_n) \longrightarrow \infty$ when $n \longrightarrow \infty.$
\end{enumerate}
Therefore $(\overline {\g}_n (k_n) -k_n)/d ( \overline {\g}_n (k_n) , k_n)  \longrightarrow w_1$ when
$n \longrightarrow \infty.$ It follows immediately from (8) that
$\alpha (\overline{\g}_n ) \longrightarrow B( v_{\theta_1(\overline {\g}_n )}, w_1) = -1.$ Therefore there exists a
$\overline {\gamma}_n$ such that
$\alpha ( \overline {\gamma}_n ) < 0.$
Let $g \in \Gamma$ be an element with
$\alpha ({g}) >0.$ If $\dim A^-(\theta_2(g)) = \dim A^+(\theta_2(\overline {\g}_n))$ set
$h = \overline{\g}^{-1}_n.$ Then $\alpha(h) < 0$ and $\dim A^-(\theta_2(g)) + \dim A^+(\theta_2(h)) =3$. Otherwise set $h = \overline{\g}_n.$ It is easy to see that there exists $t \in \Gamma$ such that
$g$ and $t h t^{-1}$ are transversal. Since $\alpha ({t h t^{-1}}) = \alpha(h)$ we have proved that there are two transversal elements in $\Gamma$ with opposite sign.\bigskip\\
\noindent\pro{ \it Proposition 6.10} Let $\Gamma$ be an affine group and let $G$ be the Zariski closure of $\Gamma.$  Assume that a semisimple part $S$ of $G$ is as in the Case 2 (3). Then $\G$ is not a crystallographic group. \endpro\\
\pro{\it Proof} \endpro This follows  immediately from Lemma 6.5 and the Main Lemma 6.9. \\
\pro{\it Remark 4.}\endpro It is possible to show that there exists an affine  group $\G  \subseteq \Aff(\mathbb{R}^6)$ acting properly discontinuously such that
the linear part of $\Gamma$ is Zariski dense in $SO(2,1) \times SL_3(\mathbb{R}).$  \bigskip\\
\pro{\it Main Theorem} Let $\G$ be a crystallographic subgroup of \text{Aff}$(\mathbb{R}^n)$ and \\ $n <7.$ Then $\G$ is virtually solvable.
\endpro\\
\noindent\pro{\it Proof .} \endpro Let $G$ be the Zariski closure of the group $\Gamma.$ Let $\dim V \leq 5.$ Then $\Gamma$ is virtually solvable by Proposition 4.2. Let $\dim V =6.$ Assume that the semisimple part $S$ of $G$ is not trivial. It follows from [S2], [To2] that the real rank of at least one simple factor group of $S$  is $\geq 2$ if $\G$ is crystallographic. Therefore $S$ is one of groups listed in Case 1 and 2. Thus $\Gamma$ is not crystallographic by Propositions 5.1, 5.2, 5.3 and 6.10. This contradiction shows that $S$ must be the trivial group. Hence the group $\Gamma$ is virtually solvable.\bigskip\\
\pro{\it Remark 5.}\endpro Actually we can prove the following proposition. Let $\Gamma$  be an affine group which acts properly discontinuously on $\mathbb{R}^n, n \leq 6$ and let $G$ be the Zariski closure of $\Gamma.$  Then $\Gamma$ is virtually solvable if and only if the complex form of any simple non-compact normal subgroup of $G$  does not have $SO_3(\mathbb{C})$ as a quotient group. 
\par Indeed,by Table 1 there is only one case left to prove; namely, $\Gamma \subseteq \text{Aff}(\mathbb {R}^6)$ and the semisimple part of the Zariski closure of $\Gamma$ is $SO(3,2).$ In this case we can prove that $\G$ does not act properly discontinuously based on dynamical ideas introduced in [AMS4.]  We show that there are two elements $g$ and $h$ of $\G$ and a compact set $K \subset \mathbb{R}^6$ such that the set $$\{n \in \mathbb{N}, m \in \mathbb{N}, k \in \mathbb{N} : g^m h^n g^m K \cap K \neq \varnothing \} $$ is infinite. Thus $\G$ does not act properly discontinuously. 
\section{The Auslander conjecture in dimension $7$. Open questions. }\label{0}
We would like to state the following  problems\\
\noindent\pro{\it Problem 1} Does there exist a crystallographic group $\G \subseteq $  \text{Aff}$(\mathbb{R}^7)$  such that $l(\G)$ is Zariski dense in $SO(4,3)$ ?  \endpro\\
We believe that this question is crucial for further progress on the Auslander conjecture. 
\par Let $G$ be the simplest representation of the simple Lie group of type $G_2$. It is well known that $G$ is a proper subgroup of  $ O(4,3).$\\
\noindent\pro{\it Problem 2} Does there exist a crystallographic group $\G \subseteq $  \text{Aff}$(\mathbb{R}^7)$  such that $l(\G)$ is Zariski dense in $G$ ?  \endpro \\
We think that these problems are very difficult. The cohomological argument used in the proof of Proposition 5.3  does not work here. This is because the virtual cohomological dimension of $ \G$ is 7 and the dimensions of the corresponding symmetric spaces are $\geq  8.$ Note that  $\alpha(\g) =\alpha(\g^{-1})$  by  6.4. Thus there is no simple way to change the sign of a hyperbolic element of $SO(4,3).$ \\
We can show that the negative answer to Problem 2 will lead to a proof of the following conjecture\\
\noindent\pro {\it Conjecture.} Let $G$ be a connected Lie group. Assume that the real rank of any simple non-commutative connected subgroup of $G$ is $\leq 2$. If a crystallographic group $\G$ is a subgroup of $G$ than $\G$ is virtually solvable.\endpro  \\
 



\begin{thebibliography}{BBNWZ}
\bibitem[A]{A} Abels, H., {\slshape Properly discontinuous groups of
affine transformations. A survey}, Geom. Dedicata 87 (2001),
309--333.
\bibitem[AMS1]{AMS1} Abels, H., Margulis  G.A., Soifer G.A.:
{\slshape Semigroups containing proximal linear maps}, Israel J.
of Math. 91 (1995), 1--30.
\bibitem[AMS2]{AMS2} ---:
{\slshape Properly discontinuous groups of affine transformations
with orthogonal linear part}, Comptes Rendus Acad. Sci. Paris 324
I (1997), 253--258.
\bibitem[AMS3]{AMS3} ---: {\slshape On the Zariski closure of the
linear part of a properly discontinuous group of affine
transformations}, Journal of Diff. Geom. 60, 2 (2003), 314-344.
\bibitem[AMS4]{AMS4} ---: {\slshape The linear part of an affine group acting
properly discontinuously and leaving a quadratic form invariant},
Geom Dedicata (2011) 153: 1-46.
\bibitem[AMS5]{AMS5} Abels, H., Margulis  G.A., Soifer G.A.:
\textit{The Auslander conjecture for dimension $ < 7$},  https://arxiv.org/pdf/1211.2525.pdf
\bibitem[Au]{Au} Auslander, L.: {\slshape The structure of complete locally affine manifolds},
Topology 3 Suppl. 1., (1964), 131-139.
\bibitem[B1]{B1} Bieberbach, L.:{\slshape \"Uber die Bewegungsgruppen des n-dimensionalen euklidischen Raumes mit einem endlichen Fundamentalbereich},
Nachr. K\"onigl. Ges. Wiss. G\"ottingen. Math. phys. Klasse (1910), 75-84.
\bibitem[B2]{B2} Bieberbach, L.:{\slshape \"Uber die Bewegungsgruppen der Euklidischen R\"aume (Erste Abhandlung)},Math. Ann. 70 (1911), 297-336. 
\bibitem[B3]{B3} Bieberbach, L.:{\slshape \" Uber die Bewegungsgruppen der Euklidischen R\"aume (Zweite Abhandlung) Die Gruppen mit einem endlichen Fundamentalbereich}, Math. Ann. 72 (1912), 400-412.
\bibitem[CDGM]{CDGM} Charette, V., Drumm, T., Goldman, W., Morrill, M. :  {\slshape Flat Affine and
Lorentzian Manifolds}, Geometriae Dedicata 97, (2003),
187 -198.
\bibitem [F] {F}E. Fedorov: {\slshape Symmetrie der regelmassigen Systeme von Figuren}, l890.
\bibitem [FK] {FK}R. Fricke and F. Klein: {\slshape Vorlesungen uber die Theoire der automorphen Functionen}, Teubner, Leipzig, 1897, Abschnitt I, Chapters 2 and 3.
\bibitem[FG]{FG} Fried, D., Goldman  W.:  {\slshape
Three--dimensional affine crystallographic groups}, Adv. in Math.
471 (1983), 1--49.
\bibitem[DG1]{DG1} Drumm, T., Goldman, W.:{\slshape Complete flat Lorenz 3--manifolds with free fundamental group}, Int. J. of Math. 1 (1990), 149-161.
\bibitem[DG2]{DG2} ---, {\slshape The geometry of crooked planes}, Topology 38 (1999), 323-351.
\bibitem[DI]{DI}  Dekimpe, K.,  Igodt, P.: { \slshape Polycyclic-by-finite groups admit a bounded-degree polynomial structure}, Invent. Math. 129 (1997), 121-140
\bibitem[GK]{GK} Goldman, W.,  Kamishima Y.: {\slshape The
fundamental group of a compact flat Lorentz space form is virtually polycyclic}, J. Differential Geom. 19 (1984), 233--240.


\bibitem[Gr]{Gr} Gromov, M.: { \slshape Almost Flat Manifolds}, J. Differential Geom., 13 (1978), 231 -241.
\bibitem[GrM]{GrM} Grunewald, F., Margulis  G.A.:
{\slshape Transitive and quasitransitive actions of affine groups
preserving a generalized Lorentz-structure}, J. Geom. Phys. 5
(1989), 493--531.
\bibitem [Hil]{Hil} Hilbert, D.: {\slshape Mathematische Probleme. Vortrag, gehalten auf dem internationalen Mathematiker-Kongre{\ss}  zu Paris 1900}. G{\"o}ttingen, 1900.
\bibitem [M]{M} Margulis, G.A.: {\slshape Complete affine locally flat manifold
with a free fundamental group}, J.Soviet Math. 134 (1987),
129--139.
\bibitem[MS]{MS}  Margulis, G.A., Soifer  G.A.: {\slshape Maximal subgroups of infinite index of linear groups}, J of Algebra, 1981, 1, 1-31
\bibitem[Me]{Me} Mess, G.: {\slshape Lorentz spacetimes of constant curvature}, Geometriae Dedicata,
126 (2007), 3-45.
\bibitem[Mi1]{Mi1} Milnor, J.: {\slshape On fundamental groups of complete affinely flat manifolds}, Adv. Math.
25 (2) (1977), 178-187.
\bibitem[Mi2]{Mi2}  Milnor, J.: {\slshape Hilbert's problem 18: On crystallographic  groups, fundamental domains, and 
on  sphere  packing}  Mathematical  Developments  Arising  from  Hilbert  Problems,”  A.M.S.  Proceedings  of  Symposia  in  Pure  Mathematics  28, 491-506,  Amer. Math. Soc.,  Providence,  R. I.,  1976. 
\bibitem[MO]{MO} Mostow, G.D.: {\slshape On the fundamental group of homogeneous spaces}, Ann. of Math. 66, 1957, 249-255.
\bibitem [OV]{OV} Onishchik, A.L., Vinberg E.B.: {\slshape Lie groups and algebraic groups}, Springer Verlag, Berlin and New York, 1990.
\bibitem[PV]{PV} Popov, V., Vinberg E.: {\slshape Invariant Theory},
Encyclopaedia of Mathematical Science 55 (1994), 123-278.
\bibitem[P]{P} Prasad, G.: {\slshape $\mathbb{R}$-regular elements in Zariski dense subgroups}, Quart. J. Math. Oxford (2) 45 (1994), 541-545 
\bibitem[R]{R}  Raghunathan, M.S.: {\slshape Discrete subgroups of Lie groups}, Springer= Verlag, Berlin and New York, 1972.
\bibitem [Ro]{Ro} K. Rohn: {\slshape Einige S{\"a}tze {\"u}ber regelm{\"a}ssige Punktgruppen}, Math. Annalen, 53 (1900), 440-449.
\bibitem [Sc]{Sc} A.Schoenflies: {\slshape Krystallsysteme und Krystallstruktur}, Teubner, Leipzig, 1891.
\bibitem[S1]{S1}  Soifer, G.: {\slshape Affine semigroup acting properly discontinuously on a hyperbolic space},
Israel Journal of Math. 192 (2012),1-20.
\bibitem[S2]{S2}  Soifer, G.: {\slshape Affine Crystallographic Groups},
 Amer. Math. Soc. Transl. (2) 163 No. 4 (1995), 165-170.
\bibitem[T]{T}  Tits, J.: {\slshape Free subgroups in linear groups},
Journal of Algebra 20 (1972), 250--270.
\bibitem[To1]{To} Tomanov, G.M.: {\slshape The fundamental group of a generalized Lorentz space form is virtually solvable }, Preprint, Tata Institute of Fundamental Research, Bombay, 1989.
\bibitem[To2]{To} Tomanov,  G.: {\slshape The virtual solvability of the fundamental group of a generalized Lorentz space form}, Journal of Diff. Geom. 32, (1990), 2, 539-547.
\bibitem[To3]{To} Tomanov, G. M.: {\slshape On a conjecture  of L. Auslander}, Comptes rendus de l'Academie bulgare des Sciences. 43, (1990), 2, 9-12.
\end{thebibliography}
\end{document}